\documentclass{article}
\usepackage[english]{babel}
\usepackage[fontsize=12pt]{fontsize}
\usepackage{amsmath}
\usepackage{graphicx}
\usepackage[colorlinks=true, allcolors=blue]{hyperref}
\usepackage{algorithm}
\usepackage{algpseudocode}
\usepackage[a4paper,top=2cm,bottom=2cm,left=2.5cm,right=2.5cm,marginparwidth=1cm]{geometry}
\usepackage{lineno,hyperref}
\usepackage{chngcntr}
\usepackage[T1]{fontenc}

\usepackage{pmboxdraw}
\usepackage{amstext}
\usepackage{amsmath}
\usepackage{amsthm}
\usepackage{graphicx}
\usepackage{epstopdf}
\usepackage{xcolor}
\usepackage{tabularx}
\usepackage{color}
\usepackage{graphicx}
\usepackage{lscape}
\usepackage{subfig}
\usepackage{subcaption}
\usepackage{float}
\usepackage{titlesec}
\usepackage{multirow}
\usepackage{chngcntr}
\usepackage{authblk}
\usepackage{caption}
\title{From Opinion Polarization to Climate Action: A Social-Climate Model of the Opinion Spectrum}
\author[1,2]{Athira Satheesh Kumar}
\author[3]{Krešimir Josić}
\author[1]{Chris T Bauch}
\author[2]{Madhur Anand}
\affil[1]{Department of Applied Mathematics, University of Waterloo, Waterloo, Ontario, Canada}
\affil[2]{School of Environmental Sciences, University of Guelph, Guelph, Ontario, Canada}
\affil[3]{Departments of Mathematics,  Biology and Biochemistry, University of Houston, Texas, USA}
\date{}

\begin{document}
\maketitle

\begin{abstract}
    
    We developed a coupled social-climate network model to understand the interaction between climate change opinion spread and the climate system and determine the role of this interaction in shaping collective actions and global temperature changes. In contrast to previous social-climate models that discretized opinions, we assumed opinions on climate change form a  continuum, and were thereby able to capture more nuanced interactions.
    The model shows that resistance to behaviour change, elevated mitigation costs, and slow response to climate events can result in a global temperature anomaly in excess of 2°C. However, this outcome could be avoided by lowering mitigation costs and increasing the rate of interactions between individuals with differing opinions (social learning). Our model is the first to demonstrate the emergence of opinion polarization in a human-environment system. We predict that polarization of opinions in a population can be extinguished, and the population will adopt mitigation practices, when the response to temperature change is sensitive, even at higher mitigation costs. It also indicates that even with polarized opinion, an average pro-mitigative opinion in the population can reduce emissions. Finally, our model underscores how frequent and unexpected social or environmental changes, such as policy changes or extreme weather events, can slow climate change mitigation. This analysis helps identify the factors that support achieving international climate goals, such as leveraging peer influence and decreasing stubbornness in individuals, reducing mitigation costs, and encouraging climate-friendly lifestyles. Our model offers a valuable new framework for exploring the integration of social and natural sciences, particularly in the domain of human behavioural change.
\end{abstract}

\section{Introduction}

Can changing human behaviour and social dynamics influence climate dynamics? To answer this question, social-climate models--models that describe both social or behavioural dynamics and climate dynamics and how the coupling of these two systems can influence each other-- are increasingly being constructed \cite{bury2019charting, moore2022determinants, shu2023determinants, beckage2022, satheesh2025climate}. Addressing the climate crisis requires closer collaboration between the social and natural sciences  \cite{moore2022determinants, shu2023determinants}. These models, as part of a broader class of coupled human-environment systems \cite{liu2007complexity,innes2013impact,henderson2016alternative,farahbakhsh2022modelling} demonstrate how such interdisciplinary approaches can pave new pathways for sustainability science.  Recent social-climate models have identified social factors such as policy, social interaction, social norms, and opinion polarization as having impacts on determining future emission trajectories \cite{perri2024socio, bury2019charting, choudhary2024climate, ahatsi2024role, kaccani2024social}. However, existing coupled social-climate models describe climate opinions as binary choices between mitigation and non-mitigation \cite{bury2019charting, menard2021, beckage2022, satheesh2025climate}. In reality, opinions lie on a spectrum \cite{castellano2009statistical, weisbuch2003interacting, favre_2024} and climate opinions, in particular, can range from strong to neutral \cite{morrison2018increasing,maibach2009global}. 

Models of opinion dynamics in human populations have been extensively studied in an uncoupled manner. These studies have elucidated the role of social influence, network interactions, and cognitive biases, in capturing the impact of social processes on opinion formation. One of the oldest continuous opinion models is the classic DeGroot model, which describes how opinions are formed in social networks \cite{degroot1974reaching}. The DeGroot model focuses on how individuals update their own beliefs by integrating the weighted opinions of their peers. Another continuous opinion model is the Friedkin-Johnsen (FJ) model, a modification of the DeGroot model that includes stubbornness \cite{friedkin1990social}, defined as resistance to changing one's opinion \cite{reynolds2023}. For instance, the FJ model was used to capture the negotiation process leading to the Paris Agreement, accounting for the stubbornness of nations during discussions \cite{bernardo2021achieving}. Continuous opinion models have been used for almost a century to describe and understand how humans interact, yet they have not been used previously in coupled social-climate models.

The weight given to the opinions of other individuals is central to opinion formation.  In opinion models, weights can be assigned to each individual depending on proximity, opinion similarity, or other characteristics. Several social experiments show that individuals usually weigh their own opinion more than the opinions of others \cite{bonaccio__2006, morin__2019, moussad__2013, david_2018}. People tend to be positively influenced by others with similar opinions \cite{ilan_yaniv_2004, kozitsin_2019, kroly_takcs__2016, moussad__2013, ivan_v__kozitsin__2020, mavrodiev_2013}. On the other hand,  experiments show that an exchange between individuals holding extremely different opinions tends to push those opinions even further apart \cite{altafini2013, bail__2018, huckfeldt__2004, nyhan__2010, myers__1970, liu__2015}. Weights can also be adaptive, meaning that they can change over time as opinions get updated, which is commonly observed in several social scenarios \cite{zahedi__2008, vicario__2016, shaw__2007, williams__2015, baronchelli__2022}. Based on these characteristics, opinion models with weights have included different features to capture the formation of opinions on a continuous spectrum \cite{mas2010individualization, ms__2013}. 

Polarization in opinion formation is frequently observed in social networks, both in theoretical models and empirical studies. Researchers have explored how exposure to opinions shared online influences opinion formation by linking it to the effects of personalization on polarization \cite{keijzer2024polarization, andersson2021dynamics}. Experiments on Facebook suggest that people adopt similar opinions and maintain distance from dissimilar opinions, which further reinforces polarization. Models are evaluated using empirical research to connect theory with real-world data \cite{flache2017models}. Factors such as politics, mass media, and social media also play an important role in shaping opinions and thus contributing to polarization \cite{wilson__2020, brulle2012shifting}. Particularly, the vast amount of information on social media introduces challenges in identifying credible sources \cite{fenoaltea__2022, fenoaltea__2023, liao2023exploring, bhattacharya2016mass, lee__2022} leading to polarization. Research also indicates that exposure to opposing opinions can increase polarization in social media \cite{bail__2018}, and the question of how polarization emerges \emph{de novo} is a topic of opinion modelling studies \cite{flache2018between,mas2014cultural,peralta2022opinion}. 

Given the widespread literature on modelling the dynamics of continuously distributed opinions, it is surprising that these approaches have not been incorporated into coupled social-climate models.  Similarly, it is surprising that the impact of environmental dynamics on continuously distributed opinion formation has not been explored in models, given that introducing human-environment coupling in binary decision models is known to create new regimes that are not possible when human or environmental systems are isolated from one another \cite{innes2013impact,sigdel2017competition,farahbakhsh2022modelling,beckage2018linking,moore2022determinants, beckage2022,bury2019charting, menard2021, satheesh2025climate}. Here, we address both knowledge gaps by modifying the FJ model with stubbornness and coupling it to an Earth system (climate) model. Stubborn individuals give a low weight to the opinions of other individuals while forming their opinions. The resulting coupled social-climate model is used to ask the following questions:  How do factors like individual stubbornness, social influence, mitigation costs, and response to temperature change alter projected emission pathways and temperature trends? How might polarized opinions emerge in a population with an initially uniform distribution of opinions? Can polarization be reduced or eradicated in a population?  How do environmental or social changes influence individual opinions and temperature trends? 

\section{Framework}
\subsection{Continuous opinion model}
We start with a population of $N$ individuals, with each individual $i$ having an opinion $o_i(t)$ about climate change. The opinion of individual $i$ in a year $t$ can take any value in $[-1,1]$, where $-1$ corresponds to a strong non-mitigative opinion, $1$ corresponds to a strong mitigative opinion and $0$ corresponds to a neutral opinion. We modify the FJ model by adding vital dynamics and making individual opinion dependent on climate dynamics, specifically through a change in the mean global temperature anomaly. 

Typically, opinion models are used to understand the formation of opinions over short periods of time. Climate change is a generational challenge; we wish to track opinion dynamics over multiple generations and hence include births and deaths when modelling the population. To do so, we randomly remove agents (death) from the population at a given rate and replace them with a new agent (birth). The opinion of the new agent is assigned by randomly choosing another individual's opinion. We also add Gaussian white noise to the evolution of opinions to account for random events such as extreme climate events, celebrity speeches, and policy changes that can influence the decision-making process. Social networks are typically sparse \cite{newman2010networks, naik2022bayesian}. Hence, we consider that an individual $i$ will only interact with $100$ other randomly selected individuals from the network of $N$ individuals at time $t$.

The continuous decision model with stubborn agents and temperature response is defined through a discrete-time formulation: 
    \begin{equation} \label{E:model}
        o_i(t+1)=\psi\left(w_{ii}o_i(t)+\lambda_i\left(\sum_{i\neq j}w_{ij}(t)o_j (t)+R(T)\right) +(1-\lambda_i) o_i(0) \right)+ \varepsilon_i(t)
    \end{equation}
where $R(T)$ is the response of opinion to temperature change (defined below), and $0<\lambda_i<1$ accounts for the susceptibility of agents to influence ($\lambda_i$ is the susceptibility of an individual to the opinions of others, and to the effects of climate change so that $1-\lambda_i$ is a measure of stubbornness) \cite{friedkin1990social}. We assume that individuals are highly influenced by similar opinions and weakly influenced by more distant ones, in accordance with empirical observations~\cite{mas2010individualization, mas2014cultural, yaniv2004receiving,kozitsin2023opinion}. Accordingly, the weights by which individual $j$ influences individual $i$ is defined as  $w_{ij} = e^{-\frac{|d_{ij}|}{A}}$  ($0\leq w_{ij}\leq1$), where $d_{ij}=|o_j(t)-o_i(t)|$, the difference in opinions between $i$ and $j$ and $A$ is the social influence of agents \cite{mas2010individualization}; smaller $A$ indicates a highly confident population that is only influenced by individuals with similar opinions, whereas a larger $A$ allows for the influence of differing opinions. The factor $\psi$ determines the time scale of opinion change to match the time scale (in years) of the climate dynamics. Since Eq~\eqref{E:model} can result in an updated opinion outside the interval $[-1,1]$, we also assume that when $o_i(t+1)<-1 \Rightarrow o_i(t+1)=-1$ and when $o_i(t+1)>1 \Rightarrow o_i(t+1)=1$. The term $\varepsilon_i(t)$ models random fluctuations in the opinion due to effects not explicitly included in the model. $R(T)$ (defined in Eq \eqref{4.4}) implicitly accounts for the outside interactions.  

\subsection{Earth System Model}
We use the following Earth system model \cite{millar2017modified}, which predicts future temperature levels depending on the atmospheric carbon dioxide levels: 
\begin{equation}
   \frac{dT_j}{dt}=\frac{q_jF-T_j}{d_j} ; T= \sum_j T_j;j=1,2 \label{2}
\end{equation}
\begin{equation}
    F=\frac{F_{\text{2x}}}{\ln(2)} \ln\left( \frac{C_0+C}{C_0} \right).
\end{equation}

Here $T$ is the temperature, $q_1$ $(0.33$ KW\textsuperscript{-1}m\textsuperscript{2}) is the thermal adjustment of the deep ocean, $q_2$ $(0.41$ KW\textsuperscript{-1}m\textsuperscript{2}) is the thermal adjustment of the upper ocean \cite{millar2017modified}, $F$ is the CO\textsubscript{2} radiative forcing, $C_0$ ($38.9$ GtCO\textsubscript{2}) \cite{masson2018,masson2019global} is the pre-industrial atmospheric CO\textsubscript{2} level, $F_\text{2x}$ $(4.5$ Wm\textsuperscript{-2}) \cite{millar2017modified} is the forcing due to CO\textsubscript{2} doubling, $d_1$ $(239$ yr) is the thermal equilibrium for the deep ocean, and $d_2$ $(4.1$ yr) is the thermal equilibrium for the upper ocean~\cite{millar2017modified}. 
\begin{figure}[H]
\centering
\includegraphics[width=15cm, height=7cm]{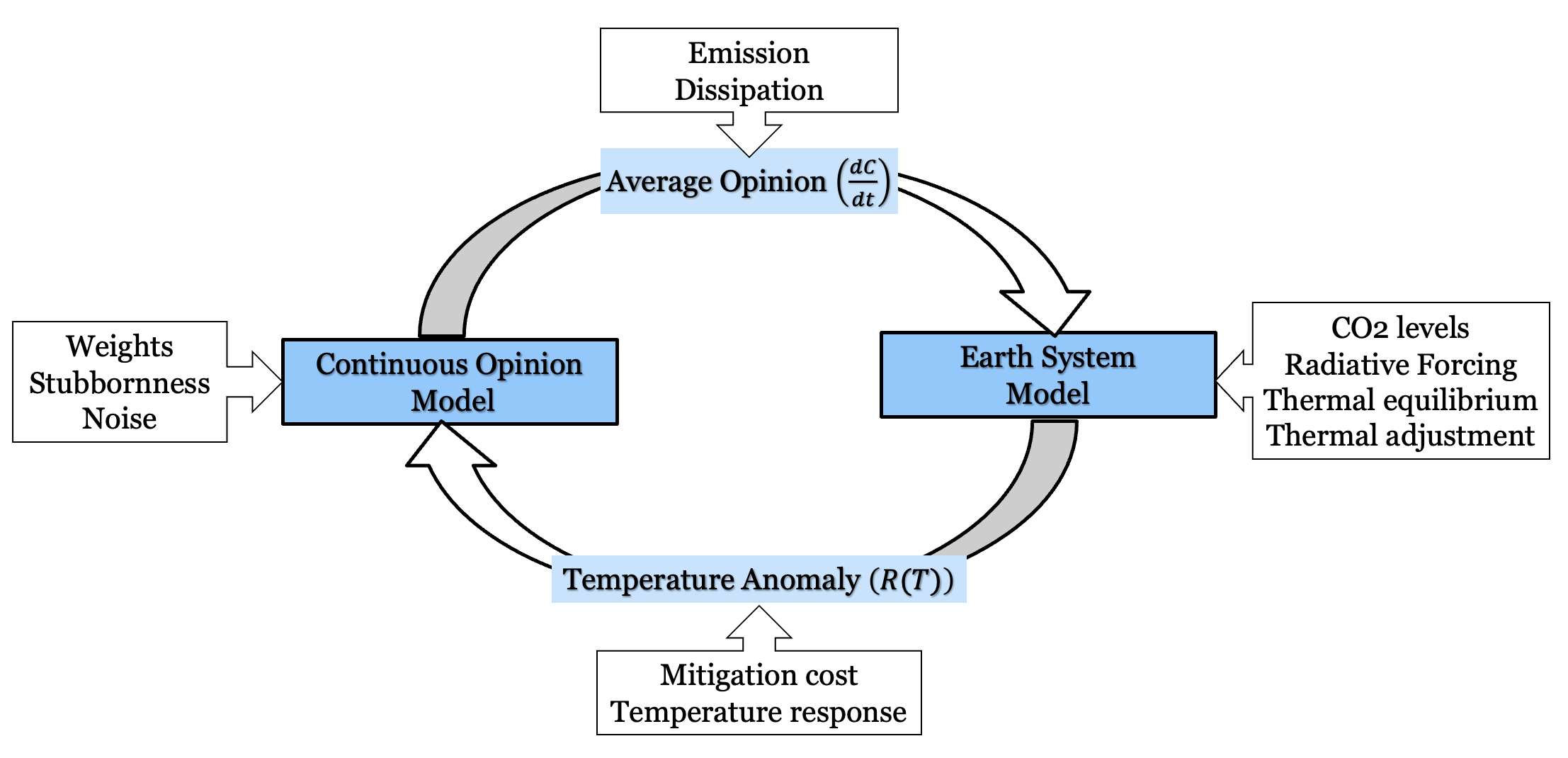}
\caption{Schematic representation of the model with input variables for each component}
\label{fig:4.6}
\end{figure}

\subsection{Coupling between opinions and climate}
We assume that all opinions,  $(-1\leq o_i \leq 1)$,  translate into actions and
that climate, in turn, impacts opinions. In particular, we model the impact of temperature $T$, on individual opinions using a sigmoid function \cite{bury2019charting, lenton2008tipping}:
    \begin{equation}
         R(T)=-m_{cost}+\left(\frac{R_{max}}{1+e^{-\alpha \left( \frac{T(t)-T_0}{T_0}\right)}}\right) \label{4.4}
    \end{equation}
depending on the cost of mitigation, $ m_{cost}> 0$, the maximum impact due to temperature change, $R_{max}$, and the sensitivity to temperature change, $\alpha$, which determines how sharply the response changes as the temperature $T$ approaches $T_0$. Here, $T_0=1.2 ^{\circ}$C is the current average global temperature anomaly against which changes are measured. We assume that every individual's emissions per unit of time are proportional to their opinion. A schematic representation of the model is shown in Fig \ref{fig:4.6}. We model the emission levels using:
\begin{equation}
    \frac{dC}{dt} =  E_\text{0} \left(0.5\left(1-\frac{\sum_{i=1}^{N} o_i(t)}{N}\right) \right)- \delta C(t) \label{4.5}
\end{equation} 
where $E_0$ is the baseline emission per year, and $\delta$ is the natural dissipation rate of carbon dioxide.   This equation captures the decrease in emissions as the average opinion moves towards acceptance of mitigation,  \emph{i.e.} as $\langle o_i(t) \rangle $ increases, which allows the opinions to influence the climate system in a stepwise manner.  We use the discrete time version of the differential equations in Eq \eqref{2} and Eq \eqref{4.5} for simulations (discrete versions in supplementary materials Eq \eqref{d1} and  Eq \eqref{d2}) to ensure consistency between the opinion and the climate model. The model chosen is discrete in time but based on its continuous version. 

\subsection{Model Parameters}

We initialize opinion dynamics in the year $2022$. 
The factor determining an individual's susceptibility to another opinion, $\lambda_i$, lies between $(0,1)$, and we assigned them by sampling values independently, uniformly at random for each agent $i$ in the population. The Gaussian white noise ($\varepsilon_{i}(t)$) added to the opinion dynamics was sampled from a normal distribution independently at each time step. For the baseline scenario, the noise is reduced to $3\%$, $\mathcal{N}(0,0.03)$, to adjust the system's sensitivity to random fluctuations.

The temperature response function $R(T)$ depends on several parameters that connect the climate system to opinion dynamics. This function depends on the cost of mitigation ($m_{cost}$), maximum response to temperature change ($R_{max}$), the sensitivity to temperature change ($\alpha$), the temperature anomaly ($T$), and the threshold temperature value ($T_0 = 1.2^\circ$ C) based on the current temperature anomaly. The values $m_{cost}$, $R_{max}$, and $\alpha$ do not have any particular range. However, these parameter values are chosen to give realistic temperature and emission predictions, which resonates with the RCP $4.5$ scenario. Representative Concentration Pathways (RCP) are various climate change scenarios that help to project the future greenhouse gas levels formally adopted by the IPCC. RCP 4.5 is an intermediate scenario where the emissions peak around 2040 and then decline \cite{meinshausen2011rcp}. By exploring various parameter combinations through trial and error, we identified that $m_{cost}=0.778$, $R_{max}=1.37$, and $\alpha = 5.7$, give peak emissions in the year 2040. The parameter that accounts for the social influence of agents, $A = 0.35$ \cite{mas2010individualization}. The scaling factor ($\psi$) determines the time scale of opinion dynamics. In the baseline scenario, we considered the scaling factor to be $0.7$. This value moderates the pace of opinion dynamics to match climate dynamics to have peak emissions occur around 2040 (consistent with the RCP 4.5 scenario), enhancing the realism of our model predictions. 

The remaining parameters have been defined previously. The pre-industrial CO\textsubscript{2} is available in the IPCC reports \cite{friedlingstein2020global}, and the natural dissipation rate is calculated from the data in \cite{friedlingstein2020global}.  We estimate the $E_0$ value from Eq \eqref{4.5}, and the data in \cite{friedlingstein2020global} by substituting all other model parameters, which gave $E_0=4.8339$ GTCO\textsubscript{2}yr\textsuperscript{-1}. 

\section{Results} \label{results}
We simulated the coupled model involving 1,000 individuals over 300 years to examine the interplay between opinion formation and climate dynamics. We treat $N=1,000$ as a representative population. Opinions are initially assigned randomly to each individual, using a truncated normal distribution, $\mathcal{N}\left(0,\frac{1}{3}\right)$ restricted to the interval $[-1,1]$ (with strong mitigative opinion represented as \(o_i=1\), strong non-mitigative opinion as \(o_i=-1\)). We then calculated the average opinion at $2150$ (at the end of our simulation) and found that it is mitigating, with a value of \(0.31\). 

\subsubsection*{A collective pro-mitigation opinion could constrain warming to $1.5^\circ$C}
The baseline parameter values result in a reduction of emission levels from \(40 \, \text{GtCO}_2 \, \text{yr}^{-1}\) to \(28 \, \text{GtCO}_2 \, \text{yr}^{-1}\) when the average opinion of the population inclines towards mitigation. This mitigating tendency in the population contributes to limiting the temperature rise beyond \(1.5^\circ\)C, which the IPCC recognizes as a critical threshold (Fig \ref{fig:4.1}(a,b,c)). The time series shows that emissions decrease slowly as soon as the opinion dynamics begin. There is a steep increase in the average opinion favouring mitigation, similar to the change in temperature levels. The opinions stabilize as the temperature starts to stabilize, which indicates the strong coupling between the average opinion and the temperature anomaly (Fig \ref{fig:B1}).

\subsubsection*{Stubborn population can cause excessive warming}
\begin{figure}[H]
\centering
\includegraphics[width=17cm, height=10cm]{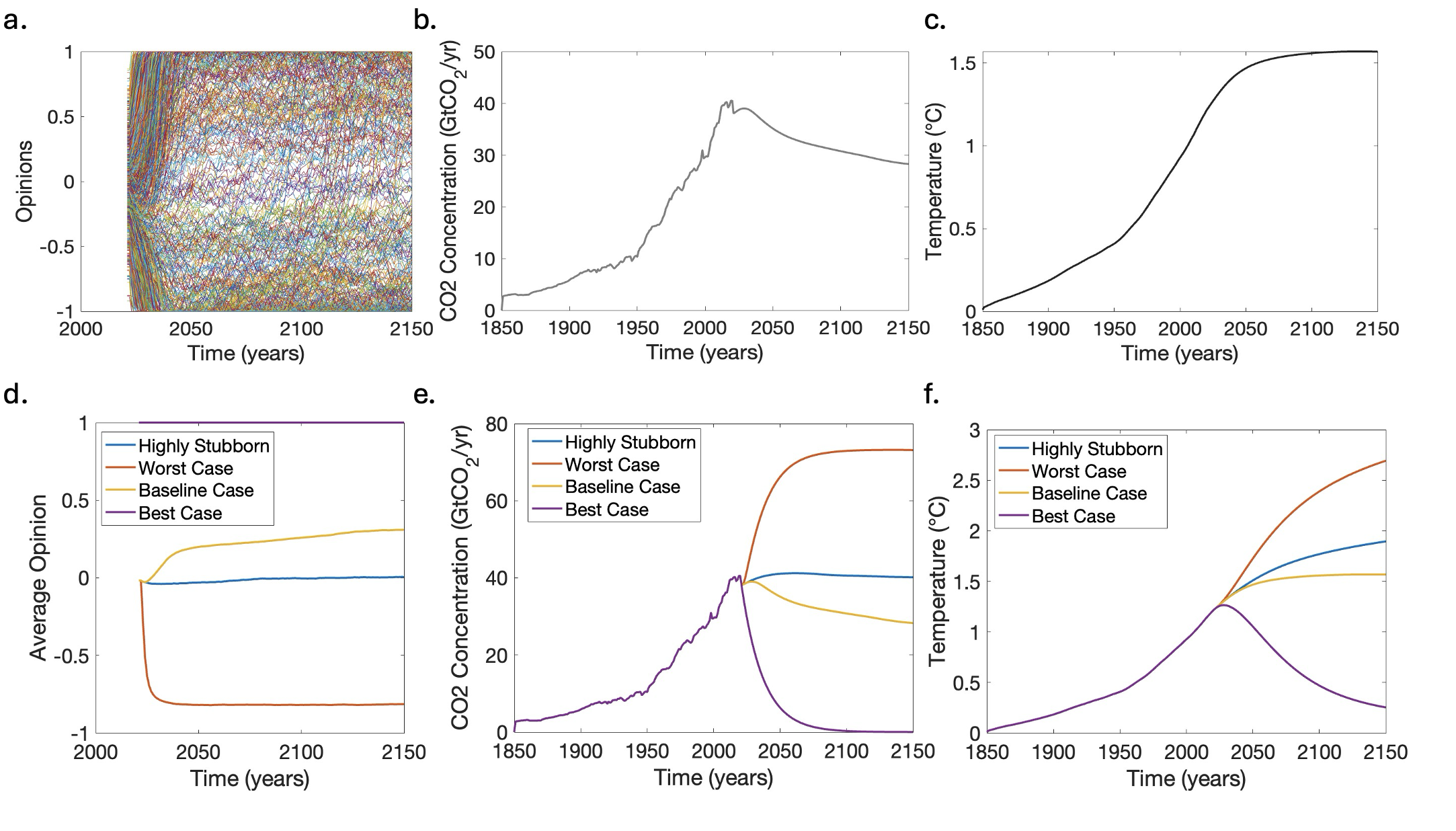}
\caption{\textbf{ Time series for the coupled social-climate model} Times series showing (a) the opinion dynamics (b) emission scenario (GtCO\textsubscript{2}yr\textsuperscript{-1}) and (c) temperature anomaly ($^\circ$C) for baseline parameter values. Time series for (d) average opinion (e) emission scenario (GtCO\textsubscript{2}yr\textsuperscript{-1}), and (f) temperature anomaly ($^\circ$C) for different parameter scenario.}
\label{fig:4.1}
\end{figure}

We next asked how stubbornness affects climate change. Specifically, suppose everyone in the population is unwilling to change their initial opinion (randomly assigned using truncated normal distribution as described above) regarding climate change. How does this shape our future climate scenario? We hypothesized that a stubborn population ($\lambda_i = 0$) reluctant to change their initial views on climate change could escalate emissions, resulting in a higher average global temperature when initial opinions are normally distributed between $[-1,1]$. We examined a scenario in which individuals are extremely stubborn ($\lambda_i = 0$, which led to the average opinion at the end of our simulation closer to $0$) and reluctant to change their perspectives (Fig \ref{fig:4.1}d). In this case, emissions remain at their current levels (Fig \ref{fig:4.1}e), while the temperature anomaly can exceed $2^\circ$C by 2100 (Fig \ref{fig:4.1}f). This suggests that a highly stubborn population unwilling to revise their decisions can trigger irreversible climate change.

We also explored the best and worst-case scenarios to understand the various possibilities of opinion, emission, and temperature scenarios. We found a best-case scenario where the entire population has strong mitigative opinions (by setting $o_i (t) =1$). In this situation, emissions drastically decline, resulting in a reduced temperature anomaly reverting to pre-industrial conditions. We further examined a worst-case scenario in which individuals refrain from choosing mitigative opinion due to a lack of response to temperature change (by setting $R_{max}=0$), resulting in the average opinion of the population moving towards an average non-mitigative opinion $(-0.8)$. Under such circumstances, emissions continue to escalate to nearly twice their existing levels, with the temperature anomaly surpassing $2.5^\circ$C by $2100$ and continuing to increase thereafter. 

\subsubsection*{Emergence of polarization}
\begin{figure}[H]
\centering
\includegraphics[width=17cm, height=10cm]{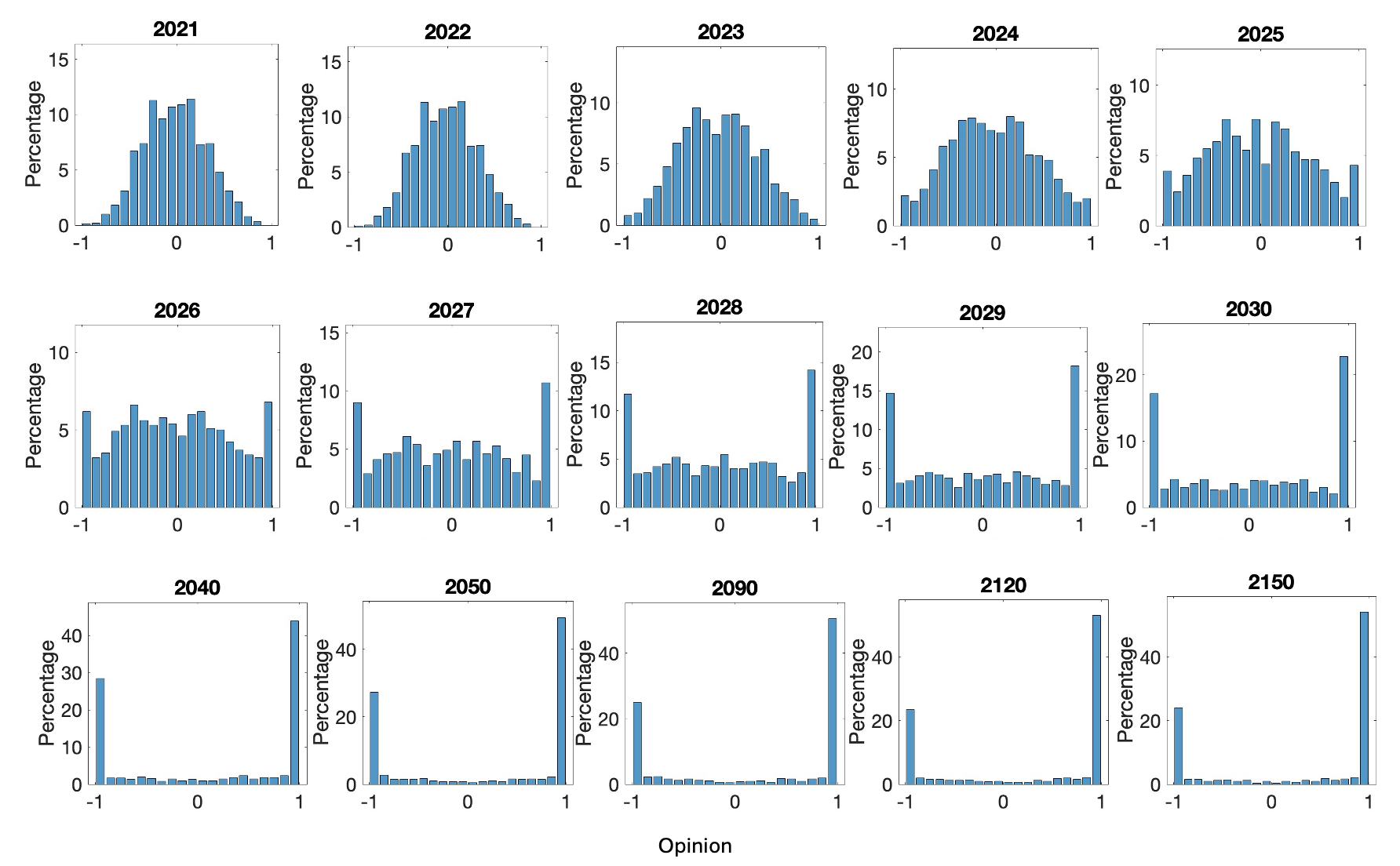}
\caption{\textbf{ Opinion evolution over time} Opinion distribution of the population at different time steps.}
\label{fig:4.2}
\end{figure}

 We did not assume that opinions are polarized \emph{a priori} but observed that polarization emerges naturally in the baseline scenario (Table \ref{table:B1}) despite the strong coupling with temperature. This observation also imply that a stronger coupling that diminishes the strength of our current temperature coupling could hinder mitigation. Although polarization emerges for a broad range of parameters in the FJ model, the coupling changes the balance in favour of mitigation. We examine the dynamics of opinion formation using the opinion distribution (Fig \ref{fig:4.2}). The initial set of opinions follows a truncated normal distribution, as mentioned above. Most opinion formation occurs within the initial decade (Fig \ref{fig:4.2}). The population is divided into highly polarized clusters of strong mitigative and non-mitigative opinions (Fig \ref{fig:4.2}). As the coupled dynamics begin, the opinion distribution gradually evolves into a flat distribution, increasing the variance of the initial normal distribution. By the year $2026$, cluster formation can be observed with less than $5\%$ of the population in each cluster. By $2027$, distinct polarization is apparent with fewer than $10\%$ of the population with strong opinions on both mitigative and non-mitigative fronts. By $2040$, the polarization intensifies with stronger mitigative opinions than non-mitigative. There are smaller subgroups between the extremities, although they constitute a small percentage($\leq 5\%$). By the end of our simulation ($2150$), more than $50\%$ of the whole population opts to adopt mitigative opinions with around $20\%$ with non-mitigative opinions and a small segment of individuals with less pronounced opinions in between.
The evolution of opinions over time is shown in Fig \ref{fig:4.2} for our chosen baseline parameters (Table \ref{table:B1}). The model dynamics indicate that polarization can occur naturally. Opinion polarization initiates after 5 years, showing slight increments as time progresses. Distinct polarization is evident 20 years after the commencement of opinion dynamics. It is also noted that an initial polarization involving equal proportions of mitigative and non-mitigative opinions changes to a state of polarization with more mitigators and fewer non-mitigators by $2100$. Our model simulation suggests that effective climate change mitigation is achievable even with $50\%$ of mitigators in the population, provided that the average opinion is mitigative. Even though extreme mitigation is not essential, reducing emissions immediately is extremely important to stop warming (Fig \ref{fig:4.1}(e,f)).

\subsection{Sensitivity Analysis}
We conducted a sensitivity analysis to understand the impact of changes in each parameter on the dynamics of the model. We executed both univariate and bivariate analyses, concentrating on how variations in model parameters influence the mean opinion, peak temperature, emission levels, and polarization. The univariate analysis is performed by changing a single parameter while maintaining all other parameters at their baseline values. This analysis helps to understand the effect of each parameter on the dynamics of the model. We investigate how these changes affect the mean opinion, emission levels, peak temperature, and polarization. In the bivariate analysis, two parameters are varied simultaneously, while all other parameters remain fixed at their baseline values. This analysis is particularly advantageous in understanding the combined effect of changes in the model parameters on the dynamics of the model. Polarization in our model is evaluated using the \emph{bimodality coefficient} (BC), which measures the degree of polarization in a distribution based on skewness and kurtosis \cite{knapp2007bimodality}. Bimodality indicates the presence of two modes in the distribution. This value provides insights into the existence of distinct subgroups or opinions. BC values that exceed 0.5 ($BC\geq 0.5$) imply considerable differences between the opinions of these groups, suggesting a notable divergence in perspectives among individuals within the model.

\subsubsection*{Temperature response and mitigation costs impact peak temperature}
\begin{figure}[H]
\centering
\includegraphics[width=18.5cm, height=10.5cm]{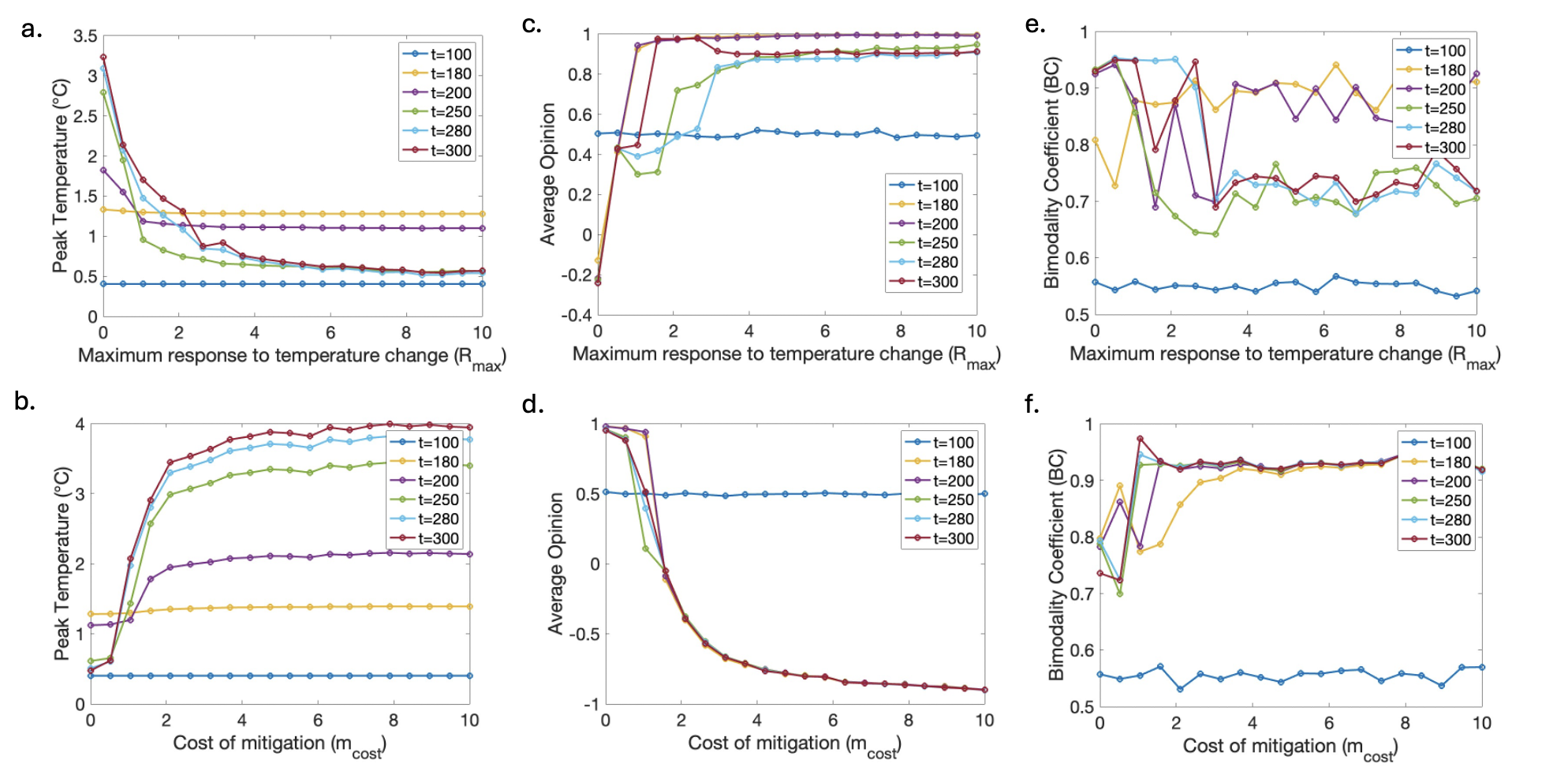}
\caption{\textbf{ Effect of variation in temperature response and cost of mitigation on opinion formation and climate change} Influence on (a) peak temperature levels, (c) average opinions, and (e) polarization using bimodality coefficient by varying $R_{max}$. Influence on (b) peak temperature levels, (d) average opinions, and (f) polarization using bimodality coefficient by varying $m_{cost}$.}
\label{fig:4.3}
\end{figure}

Univariate analysis of the parameter that determines the impact of temperature fluctuations on opinions ($R_{max}$) reveals substantial effects on climate (peak temperature, emissions) and social (average opinion and polarization) dynamics (Fig \ref{fig:4.3} (a,c,e)). A heightened sensitivity to peak temperature and average public opinion is noted when the response to temperature changes is minimal ($R_{max} \leq 2$). Decreasing $R_{max}$ can increase peak temperature ranging from $0.7$ to $3.2$ degrees within the projected future scenarios. Furthermore, a gradual transition in the average opinion of the population from a non-mitigative opinion to one advocating for mitigation is observed with an intensified response to temperature change. Our model further suggests that elevated response values enhance mitigative opinions, reducing emissions and peak temperature levels, as anticipated.
A similar examination was conducted by altering the cost of mitigation ($m_{cost}$) (Fig \ref{fig:4.3} (b,d,f)). An increase in mitigation costs is associated with a notable increase in peak temperature values (ranging from $0.5$ to $4$ degrees) with a shift in the general opinion towards non-mitigative approaches. The elevated costs associated with mitigation reduce mitigation efforts, resulting in more individuals choosing non-mitigative strategies resulting in an increase in emissions levels. Furthermore, a reduction in polarization (measured by a decreased bimodality coefficient) is observed in contexts where mitigation costs are relatively low. In both scenarios with $R_{max}$ and $m_{cost}$, more changes are evident within the parameter values spanning from $0$ to $2$ (Fig \ref{fig:4.3}).

\subsubsection*{Increased perturbations delay mitigation}
Incorporating random noise ($\varepsilon_i (t)$) into the framework of opinion dynamics introduces stochasticity that influences the evolution of individual beliefs over time. Random noise represents environmental (extreme weather events like heatwaves and floods) or social (policies or government) changes that can influence individual opinions. We hypothesized that individuals are more likely to reassess and alter their opinions when these perturbations are larger. To test this hypothesis, we conducted a univariate analysis by increasing the intensity of random noise ($\varepsilon_i(t)$). As the intensity of the random noise increases, the peak temperature shows a corresponding rise, while the mean opinion becomes less mitigative. Specifically, random noise alone was sufficient to elevate the peak temperature by $0.6^\circ$C, as illustrated in Fig \ref{fig:B4}. This indicates that regular unexpected events (environmental, political, economic, and others) obstruct mitigation efforts as they can create uncertainty in long-term planning, ultimately delaying emission reduction. This finding highlights the influence that random factors can exert on the overall temperature behaviour within the system, even though these random events might seem neutral.

\subsubsection*{Social learning considerably reduces peak temperature anomaly}
We analyzed the scaling parameter ($\psi$) that influences stubbornness and susceptibility to other opinions in the network by varying it from $0.1$ to $1$. This scaling factor is important in shaping the average opinion and temperature anomaly. We found that an increase in the scaling factor results in a notable decline in peak temperature, dropping from $2.8^\circ$C to $1.5^\circ$C. At the same time, this modification shifted the mean opinion from a neutral position towards mitigation, as illustrated in Fig \ref{fig:B2}. This interaction among scaling, opinion dynamics, and temperature signifies the connections that govern social agreement and environmental reactions to climate change.
\begin{figure}[H]
\centering
\includegraphics[width=17cm, height=7cm]{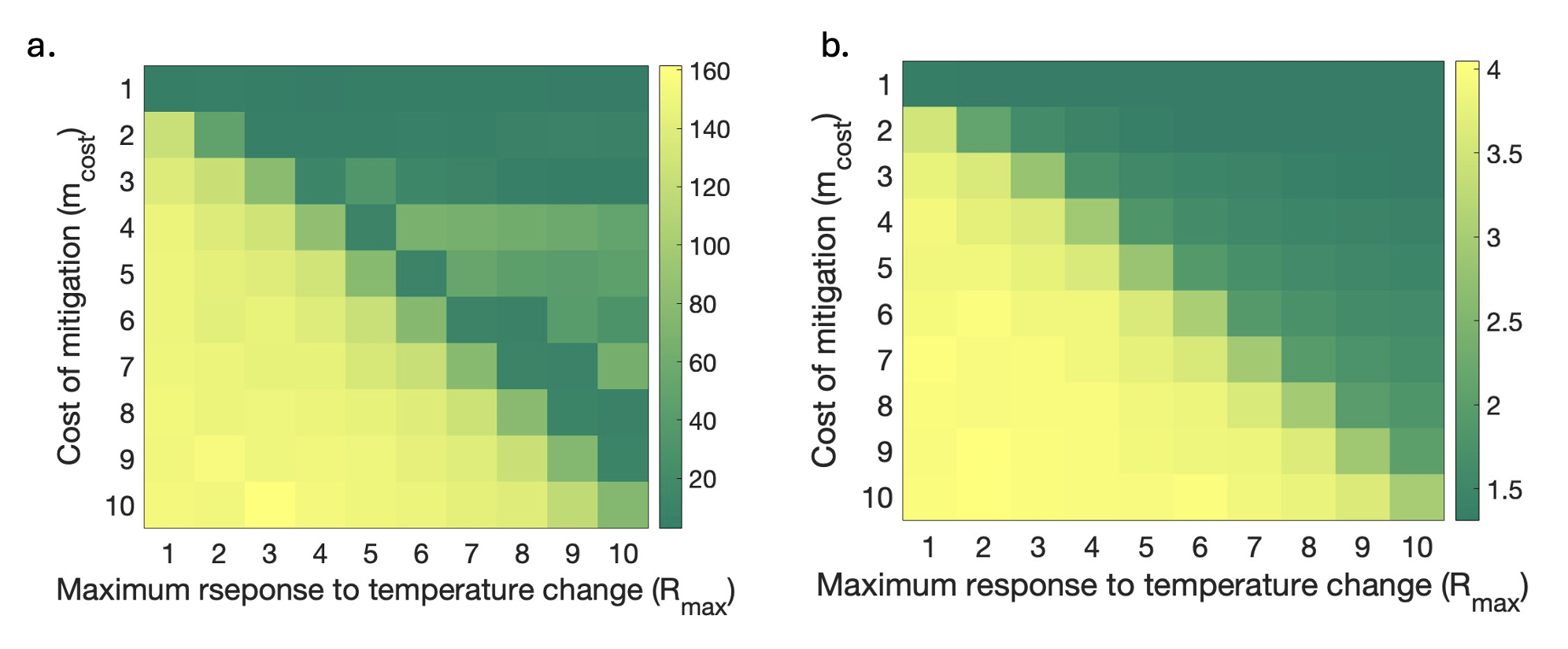}
\caption{\textbf{ Temperature response and cost of mitigation have a major influence on emission and peak temperature} Heatmap showing the influence of variation in $R_{max}$ and $m_{cost}$ on (a) emission levels and (b) peak temperature anomaly.}
\label{fig:4.4}
\end{figure}
\subsubsection*{High mitigation costs with weak temperature response increases polarization}
Conducting a bivariate analysis enabled us to evaluate and differentiate the impacts of changes in parameters on various important outputs like the peak temperature, emission rates, average opinion, and the level of polarization. Our bivariate analysis involving the variations in the parameters $R_{max}$ and $m_{cost}$, revealed that both parameters had a similar influence on the outcomes we studied, illustrating that their influences are comparable in significance. We observed that a combination of relatively low costs with higher temperature responses led to lower peak temperatures and emissions, as shown in Fig \ref{fig:4.4}, while simultaneously fostering mitigation through higher average opinions and an increased degree of polarization, as illustrated in Fig \ref{fig:4.5}. This indicates that lowering the cost of mitigation and increasing the response to temperature can effectively diminish emissions and polarization.

\subsubsection*{Heightened response to temperature changes drives unified mitigation efforts}
\begin{figure}[H]
\centering
\includegraphics[width=17cm, height=7cm]{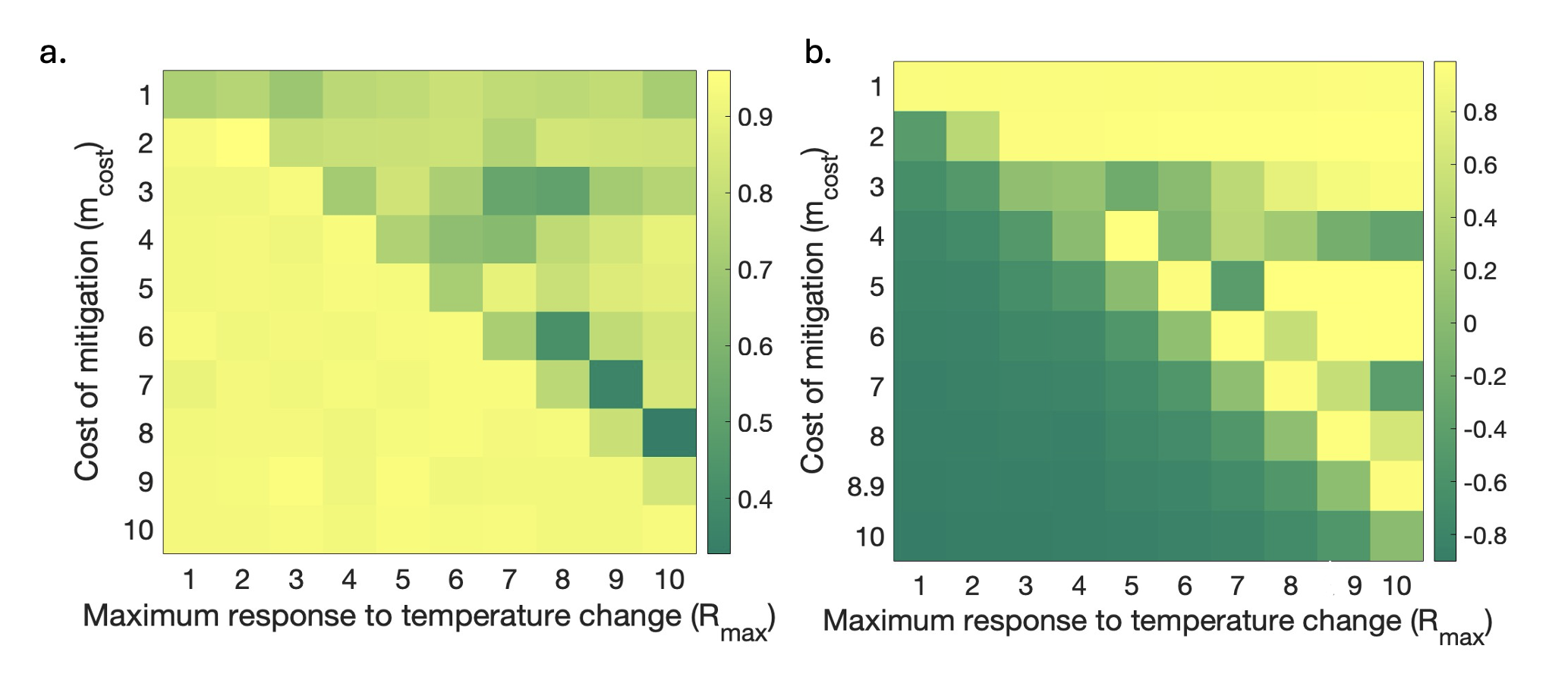}
\caption{\textbf{ Temperature response and cost of mitigation have a major influence on opinion formation} Heatmap showing the influence of variation in $R_{max}$ and $m_{cost}$ on (a) polarization using the bimodality coefficient and (b) average opinion.}
\label{fig:4.5}
\end{figure}

Comparing Fig \ref{fig:4.4}a and \ref{fig:4.5}a shows a crucial pattern relating emission levels with polarization. Despite the general trend revealing that high mitigation costs with low response increase peak temperature levels, emissions, polarization and decrease average opinion, the emission levels and polarization act differently at certain specific points. Low polarization is observed with an increased response to temperature change even when the mitigation cost is high, which suggests the possibility of unified mitigation when the individual response to mitigate these changes. The high cost of mitigation is no longer a concern for individuals at this stage. This indicates a threshold effect in which a specific interplay of $R_{max}$ and $m_{cost}$ results in reduced polarization and temperatures.

\subsubsection*{$R_{max}$ and $m_{cost}$ considerably influence model dynamics}
A bivariate analysis conducted to assess the influence of the scaling factor ($\psi$) along the cost of mitigation ($m_{cost}$), shows that the effect of the scaling factor is considerably reduced when compared to the overall impact the cost of mitigation has on the peak temperature anomaly. In particular, in circumstances where the cost of mitigation remains relatively low ($m_{cost} \leq 2$), the scaling factor has a considerable influence over the average opinion of the population; however, as the cost of mitigation begins to increase, the scaling factor’s impact on the average opinion of the population is reduced (Fig \ref{fig:B7} (a,b) and \ref{fig:B8} (b)). In this model, the scaling factor is linked to the social learning component typically incorporated in social-climate models that aim to understand human behaviour and climate change. An increase in social learning enhances efforts toward mitigation, similar to the findings in various discrete social-climate models \cite{bury2019charting, menard2021, beckage2018linking}. 

An examination focusing on the connection between the scaling factor and the temperature change response, denoted as ($R_{max}$), revealed that heightened social learning, as signified by the scaling factor, can elevate peak temperatures when the response is constrained ($R_{max} \leq 2$) (Fig \ref{fig:B5}). As social learning rises with reduced temperature responses, the collective opinion of the population transitions from an opposition to mitigation to one that supports mitigation (Fig \ref{fig:B6} (b)). In essence, the findings indicate that the scaling factor exerts a stronger influence on peak temperatures, emission rates, and the opinion of the population when both mitigation costs and temperature responses are small. An examination of polarization while adjusting these parameters shows that polarization is directly tied to the scaling factor rather than being a consequence of variations in the temperature response (Fig \ref{fig:B6} (a)). Conversely, the scaling factor appears to exert minimal influence when the cost of mitigation is higher (Fig \ref{fig:B8} (a)). Furthermore, the analysis of the impact of additional model parameters with fluctuations in $R_{max}$ and $m_{cost}$ suggests that elements such as the social influence of agents ($A$) and the sensitivity to temperature variations $(\alpha)$ have a considerably smaller effect on the overall dynamics compared to $R_{max}$ and $m_{cost}$ (Fig \ref{fig:B9} - \ref{fig:B16}). Our findings suggest that enhancing the response to temperature change and lowering the cost of mitigation is ideal for addressing climate change.

\section{Discussion}
Our model highlights the role of opinion dynamics on social networks and coupled social-climate feedbacks in shaping views on climate change and thus intensifying or mitigating climate change impacts. Climate opinions reflect not just personal beliefs but also the influence of external factors such as social interaction, individual responses to climate events, and the costs associated with mitigation. Our model shows how opinion polarization can arise in a broad parameter regime for all three of these factors. Moreover, opinion formation depends on the climate scenario and individual interactions. 

This is the first demonstration of the emergence of opinion polarization in a human-environment system, to our knowledge. Our simulations also indicate that even with polarized opinions, an average pro-mitigation opinion in the population can enhance mitigative behaviour and reduce emissions. Although having an entire population with a strongly pro-mitigation opinion is not essential, a minimum of $50$ percent of the population embracing mitigation is necessary to limit our emissions and halt warming. The ideal way to boost the number of mitigators in the population is to lower the mitigation costs and increase individual response to climate change by adopting mitigation-friendly lifestyles. However, our simulation also points out some critical values of mitigation costs that could eliminate polarization. These critical points suggest that a heightened response to climate change can lessen the concerns over the cost of mitigation at certain points. 

Feedback from the environmental systems influences long-term opinion distribution. According to our analysis of the climate change issue, frequent and unexpected social or environmental changes could also result in lower average pro-mitigation opinion in the population, leading to elevated peak temperatures.  Our findings show that social learning (represented by the scaling parameter) is important in forming opinions on climate change, and that a high social learning rate can reduce temperature anomalies. However, a highly stubborn population can result in an average global temperature anomaly in excess of $2 ^{\circ}$C. Elements such as increased responsiveness of individuals to climate change, reduced mitigation costs, and enhanced social learning can promote mitigation, especially in tandem.

Although this model offers a new perspective on coupled social-climate dynamics, it does not include all relevant aspects. Several social factors might contribute to opinion shifts regarding climate change, including climate education, media representation, economic circumstances, and political influence \cite{wilson__2020, brulle2012shifting, bernardo2021achieving}. These factors were not explicitly represented in our social learning model. There is a possibility for global consensus on certain mitigative strategies, free from biases related to politics, geography, or cultural influences \cite{bernardo2021achieving}. This could give rise to other issues when individuals do not express their true opinions to align with the prevailing political, cultural, or geographical ideologies. Incorporating these elements could serve as a valuable extension to the model, enabling deeper exploration of several aspects of society that may promote or hinder mitigation. Even without these aspects, simple coupled models like the one presented in this work can provide valuable insights into the interplay between opinion and climate systems.  More sophisticated social-climate models that build on these insights could help policymakers create targeted strategies that address prevailing opinions, ultimately enhancing the effectiveness of climate change mitigation initiatives.

\bibliographystyle{unsrt}
\bibliography{sample}
\newpage
\appendix
\section{Appendix}
\subsection{Table of Parameters}
\begin{table}[h]
 \footnotesize
\centering
\begin{tabular}{|c|c|c|c|c|} 
 \hline
Parameters& Description & Values&Units&Sources \tabularnewline
\hline 
\hline 
$N$ &  Total number of individuals &1000& 1 & calibrated\tabularnewline
\hline 
$w_{ij}$ &  Weight given to other opinions &[0,1]& 1 & \cite{mas2010individualization}\tabularnewline
\hline 
$\lambda_i$ & Susceptibility to other opinions &(0,1)& 1 & \cite{friedkin1990social} \tabularnewline
\hline 
$o_i$ & Opinion of individuals &[-1,1]&1 &\cite{mas2010individualization}
 \tabularnewline
\hline
$C_0$ & Pre-industrial $CO_2$ level &38.9& GTCO\textsubscript{2}&\cite{friedlingstein2020global}
 \tabularnewline
\hline
$F_{2x}$ & Forcing due to $CO_2$ doubling &4.5& Wm\textsuperscript{-2} &\cite{masson2021ipcc}
 \tabularnewline
\hline
$q_1$ & Thermal adjustment of deep ocean &0.33& KW\textsuperscript{-1}m\textsuperscript{2} 
&\cite{millar2017modified}
 \tabularnewline
\hline
$q_2$ & Thermal adjustment of upper ocean &0.41& KW\textsuperscript{-1}m\textsuperscript{2} &\cite{millar2017modified}
 \tabularnewline
 \hline
$d_1$ & Thermal equilibrium for deep ocean &239& yr &\cite{millar2017modified}
 \tabularnewline
\hline
$d_2$ & Thermal equilibrium for upper ocean &4.1& yr &\cite{millar2017modified}
 \tabularnewline
\hline
$m_{cost}$ & Cost of mitigation &0.778&1& calibrated\tabularnewline
\hline 
$R_{max}$ & Maximum response to temperature change &1.37& 1 &calibrated\tabularnewline
\hline
$\delta$ & Natural dissipation rate &0.06& yr\textsuperscript{-1} &\cite{friedlingstein2020global}\tabularnewline
\hline
$A$ & Social influence of agents &0.35& 1&calibrated\tabularnewline
\hline
$\alpha$ & Sensitivity to change in temperature  &5.7& 1 & calibrated \tabularnewline
\hline
$\psi$ & Scaling for susceptibility and stubbornness &0.7& 1&calibrated \tabularnewline
\hline
\end{tabular}
\caption{Table describing the model parameters and baseline values}
\label{table:B1}
\end{table}
\section{Discrete time equations for ESM}

We use the equivalent discrete time equations for the ODEs used in the ESM and coupled model to ensure consistency between the opinion and climate models.
\begin{equation}
   T_j(t+1)=T_j(t)+ \frac{q_jF-T_j (t)}{d_j} ; T= \sum_j T_j;j=1,2 \label{d1} \tag{S1}
\end{equation}
is the discrete time form of Eq \eqref{2} and
\begin{equation}
    C(t+1) = C(t)+ 0.5 E_\text{0} \left(1-\frac{\sum_{i=1}^{N} o_i(t)}{N}\right)- \delta C(t) \label{d2} \tag{S2}
\end{equation} 
is the equivalent discrete time form of Eq \eqref{4.5} used for simulation.

\newpage
\setcounter{figure}{0}
\counterwithin{figure}{section}
\subsection{Additional Figures}
\begin{figure}[H]
\centering
  \includegraphics[width=17cm, height=12cm]{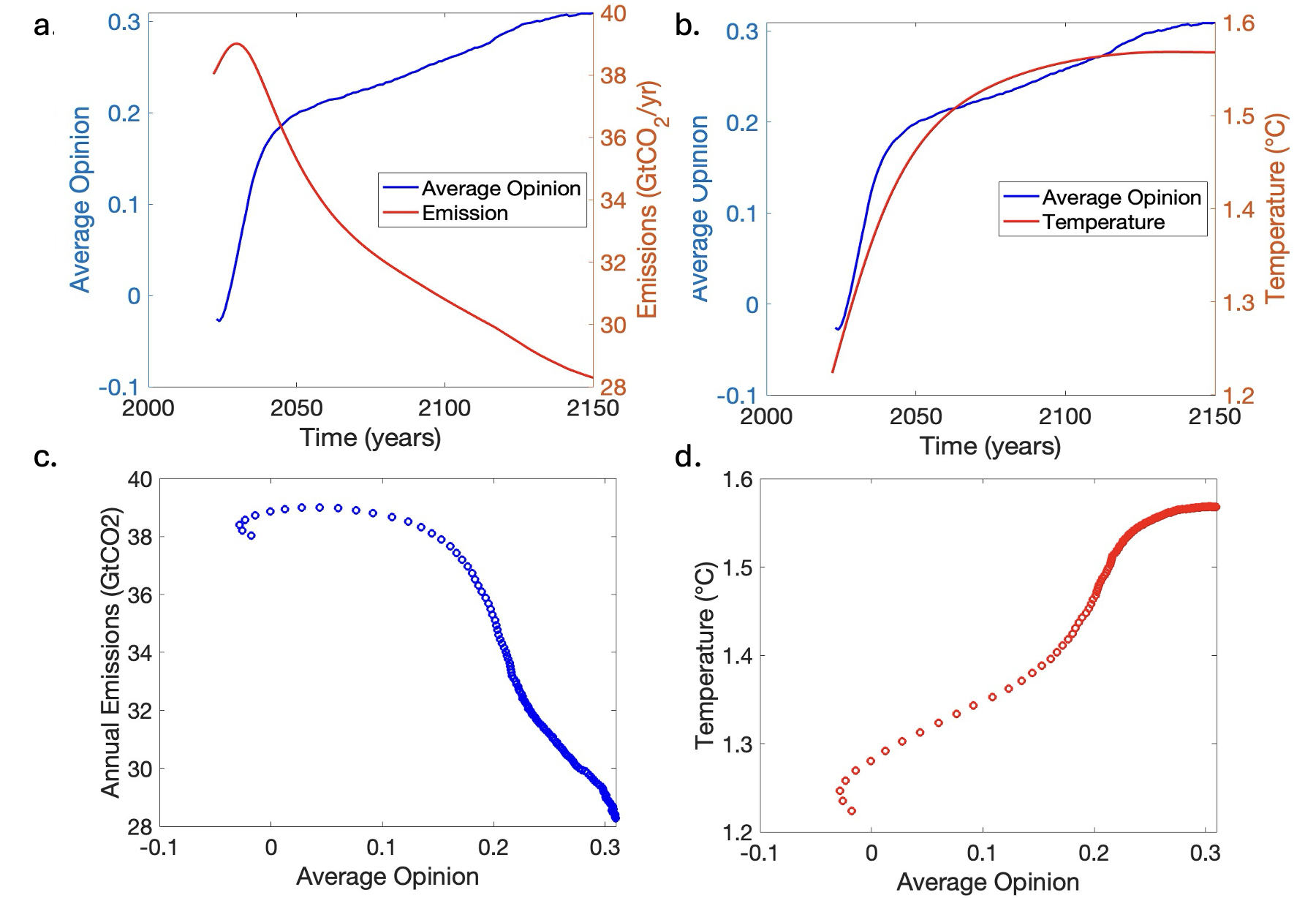}
\caption{\textbf{Relation between average opinion, emission, and temperature:}Time series of (a) average opinions and emissions (GTCO\textsubscript{2}yr\textsuperscript{-1}), (b)average opinion and temperature ($^\circ$C), and correlation between (c) average opinions and emissions (GTCO\textsubscript{2}yr\textsuperscript{-1}), (d) )average opinion and temperature ($^\circ$C).}
\label{fig:B1}
\end{figure}
\begin{figure}[H]
\centering
  \includegraphics[width=18cm, height=7cm]{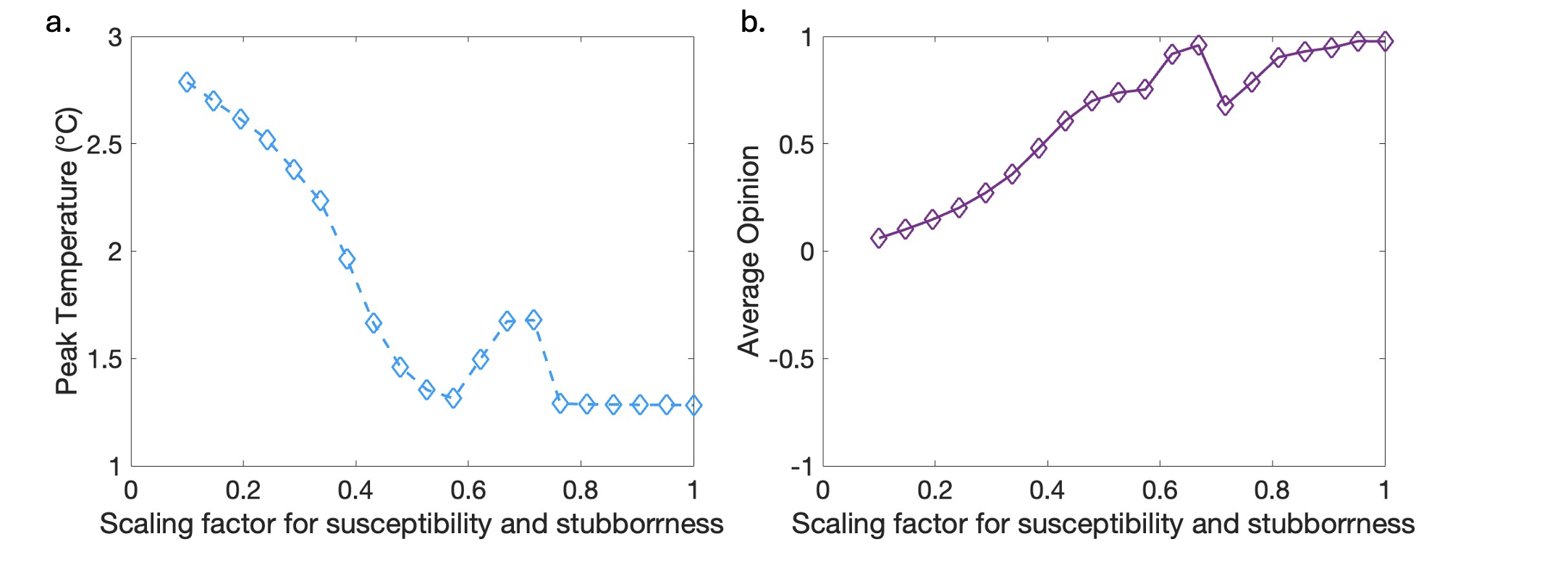}
\caption{\textbf{ High social learning mitigates climate change: }Univariate sensitivity of (a) peak temperature ($^\circ$C), (b) average opinion with variation in scaling factor of susceptibility and stubbornness.}
\label{fig:B2}
\end{figure}
\begin{figure}[H]
\centering
  \includegraphics[width=18cm, height=7cm]{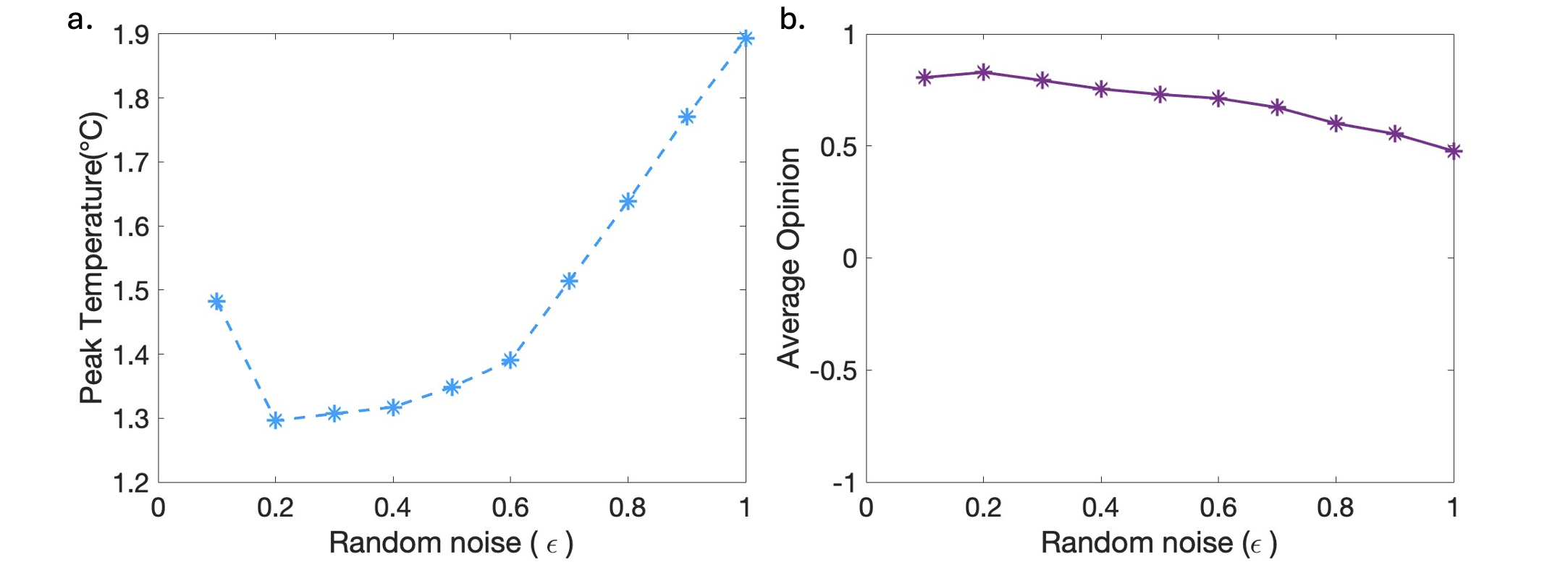}
\caption{\textbf{Increased stochasticity increases the peak temperature: }Univariate sensitivity of (a) peak temperature ($^\circ$C), (b) average opinion with variation in random noise added to the opinion dynamics.}
\label{fig:B4}
\end{figure}
\begin{figure}[H]
\centering
  \includegraphics[width=18cm, height=7cm]{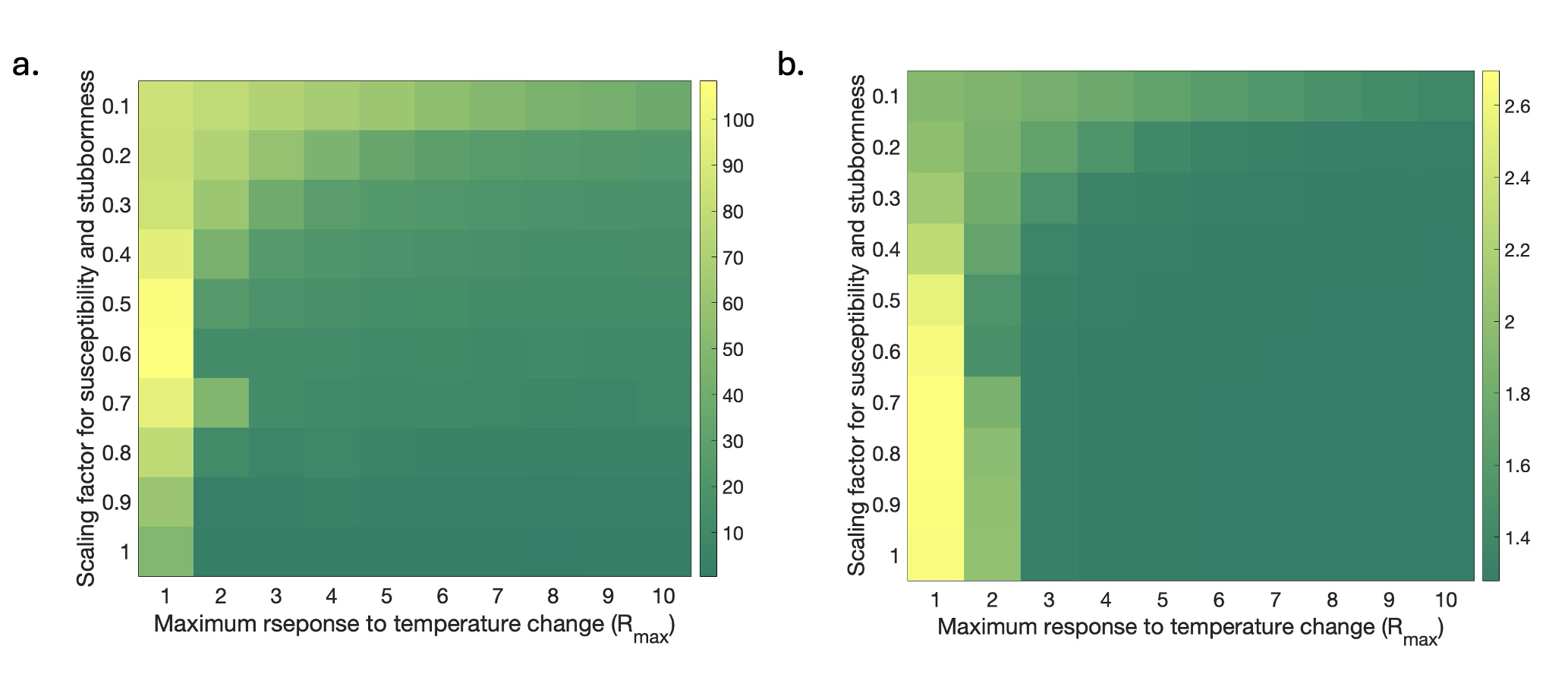}
\caption{\textbf{Scaling the opinion has minimal influence with high response: }Bivariate sensitivity of (a) emissions (GTCO\textsubscript{2}yr\textsuperscript{-1}), (b) peak temperature ($^{\circ}$C) when scaling factor is varied with the maximum response to temperature change $(R_{max})$}
\label{fig:B5}
\end{figure}
\begin{figure}[H]
\centering
  \includegraphics[width=18cm, height=7cm]{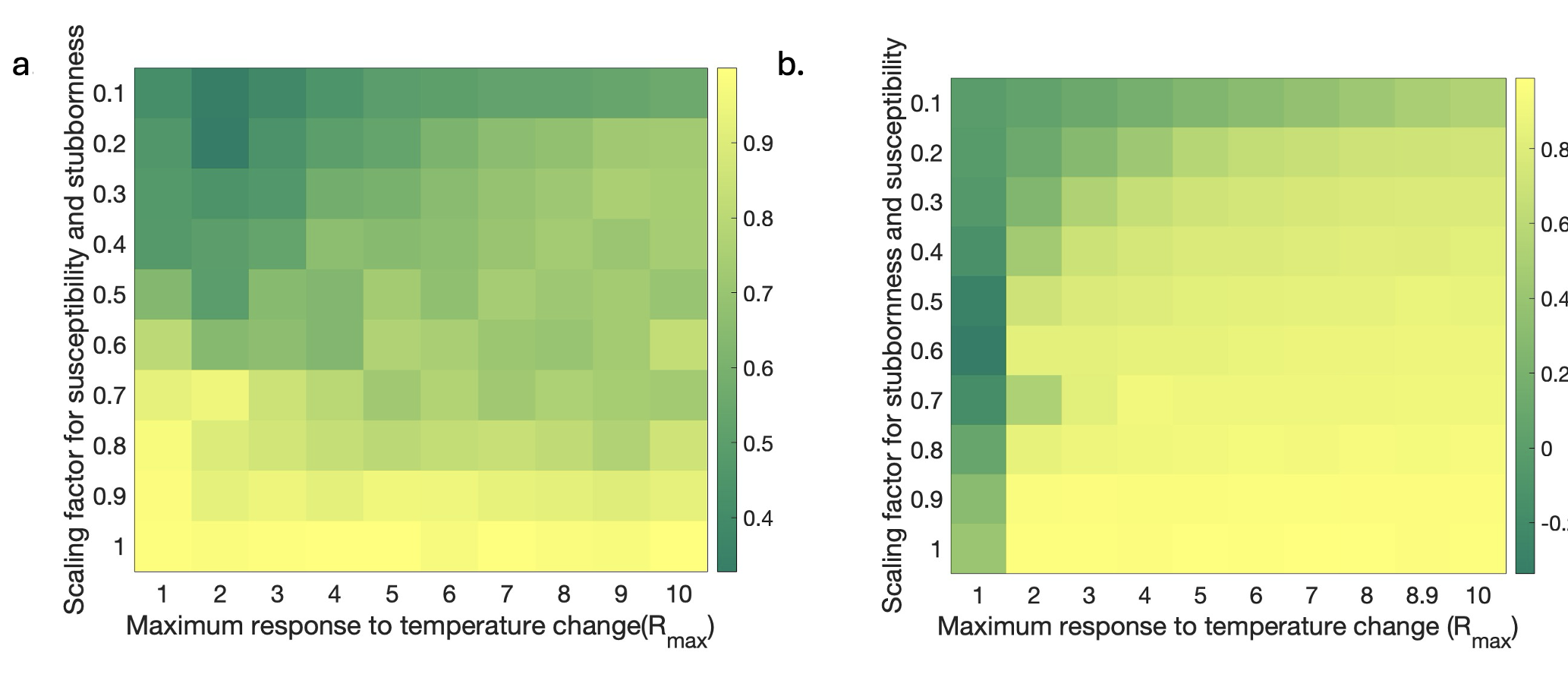}
\caption{\textbf{Scaling the opinion influences polarization: }Bivariate sensitivity of (a) bimodality coefficient, (b) average opinion when scaling factor is varied with the maximum response to temperature change $(R_{max})$}
\label{fig:B6}
\end{figure}
\begin{figure}[H]
\centering
  \includegraphics[width=18cm, height=7cm]{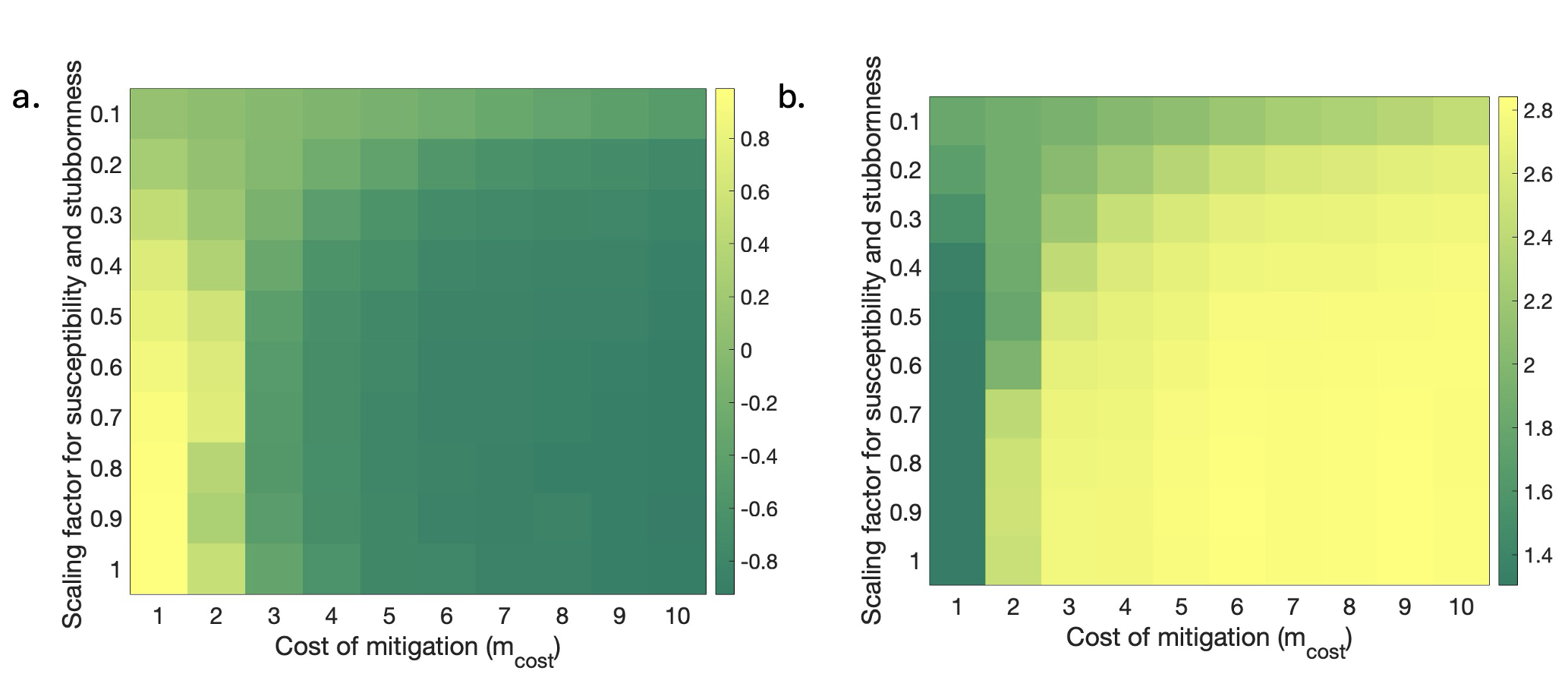}
\caption{\textbf{Scaling the opinion has minimal influence with high mitigation cost: }Bivariate sensitivity of (a) emissions (GTCO\textsubscript{2}yr\textsuperscript{-1}), (b) peak temperature ($^{\circ}$C) when scaling factor is varied with the cost of mitigation $(m_{cost})$}
\label{fig:B7}
\end{figure}
\begin{figure}[H]
\centering
  \includegraphics[width=18cm, height=7cm]{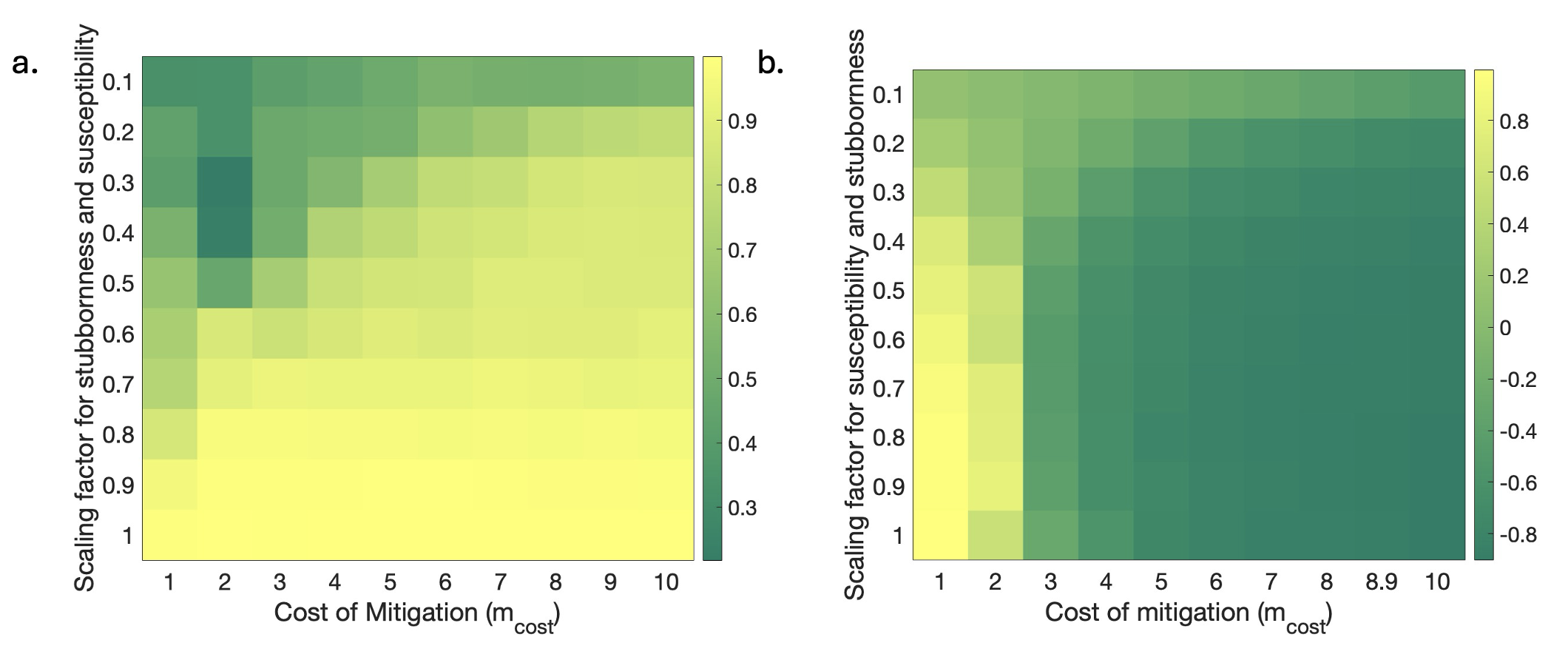}
\caption{\textbf{Scaling the opinion influences polarization with mitigation cost: }Bivariate sensitivity of (a) bimodality coefficient, (b) average opinion when scaling factor is varied with the cost of mitigation $(m_{cost})$}
\label{fig:B8}
\end{figure}
\begin{figure}[H]
\centering
  \includegraphics[width=18cm, height=7cm]{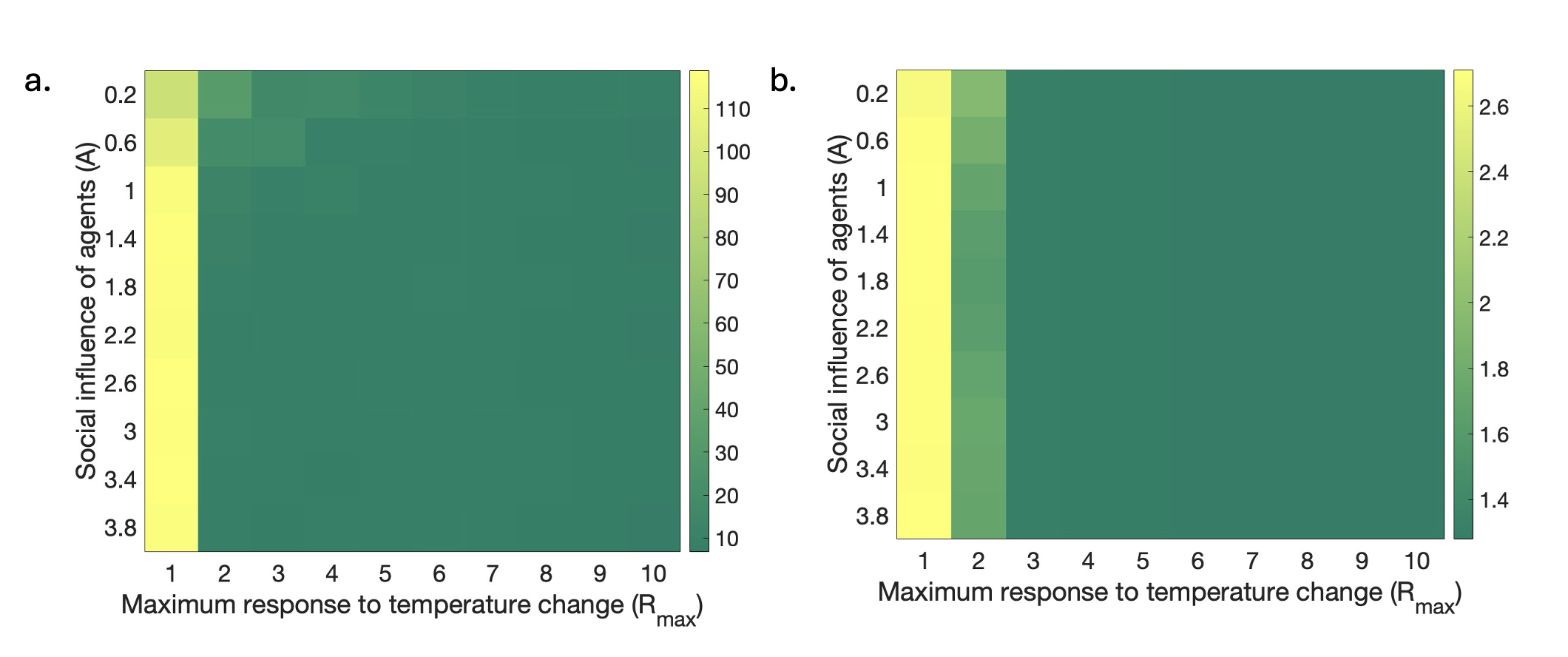}
\caption{\textbf{Social influence has minimal influence on emissions with high response: }Bivariate sensitivity of (a) emissions (GTCO\textsubscript{2}yr\textsuperscript{-1}), (b) peak temperature ($^{\circ}$C) when social influence($A$) is varied with the maximum response to temperature change $(R_{max})$}
\label{fig:B9}
\end{figure}
\begin{figure}[H]
\centering
  \includegraphics[width=18cm, height=7cm]{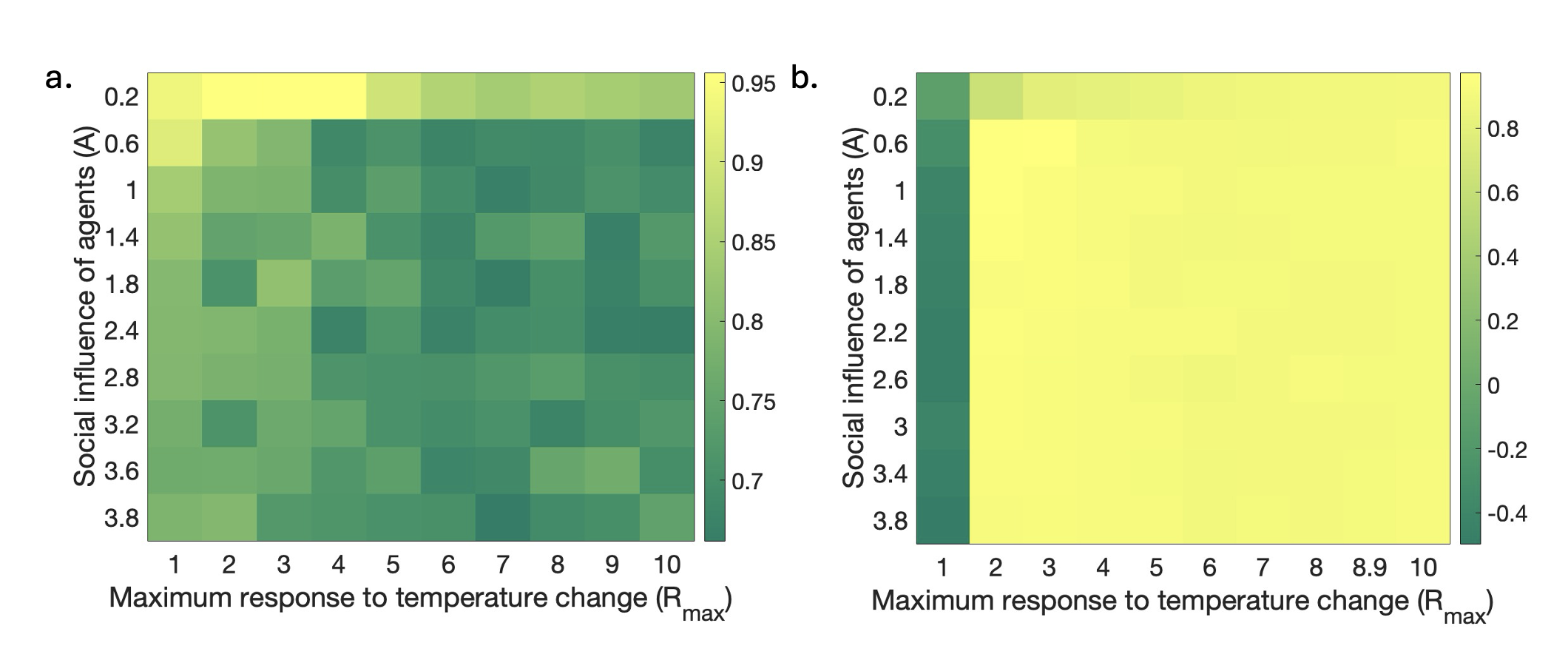}
\caption{\textbf{Social influence has minimal influence on average opinion with high response: }Bivariate sensitivity of (a) bimodality coefficient, (b) average opinion when social influence($A$) is varied with the maximum response to temperature change $(R_{max})$}
\label{fig:B10}
\end{figure}
\begin{figure}[H]
\centering
  \includegraphics[width=18cm, height=7cm]{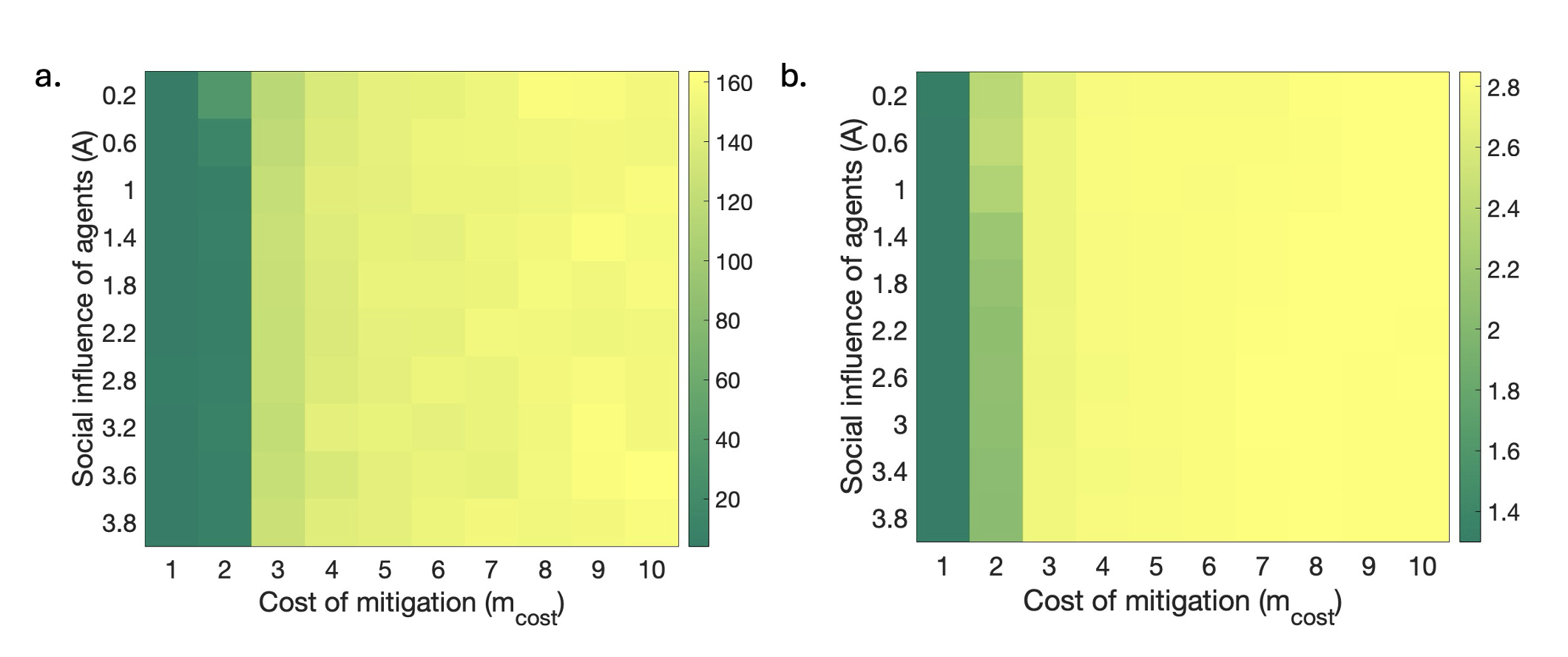}
\caption{\textbf{Social influence has minimal influence on emissions with mitigation costs: }Bivariate sensitivity of (a) emissions (GTCO\textsubscript{2}yr\textsuperscript{-1}), (b) peak temperature ($^{\circ}$C) when social influence($A$) is varied with the cost of mitigation $(m_{cost})$}
\label{fig:B11}
\end{figure}
\begin{figure}[H]
\centering
  \includegraphics[width=18cm, height=7cm]{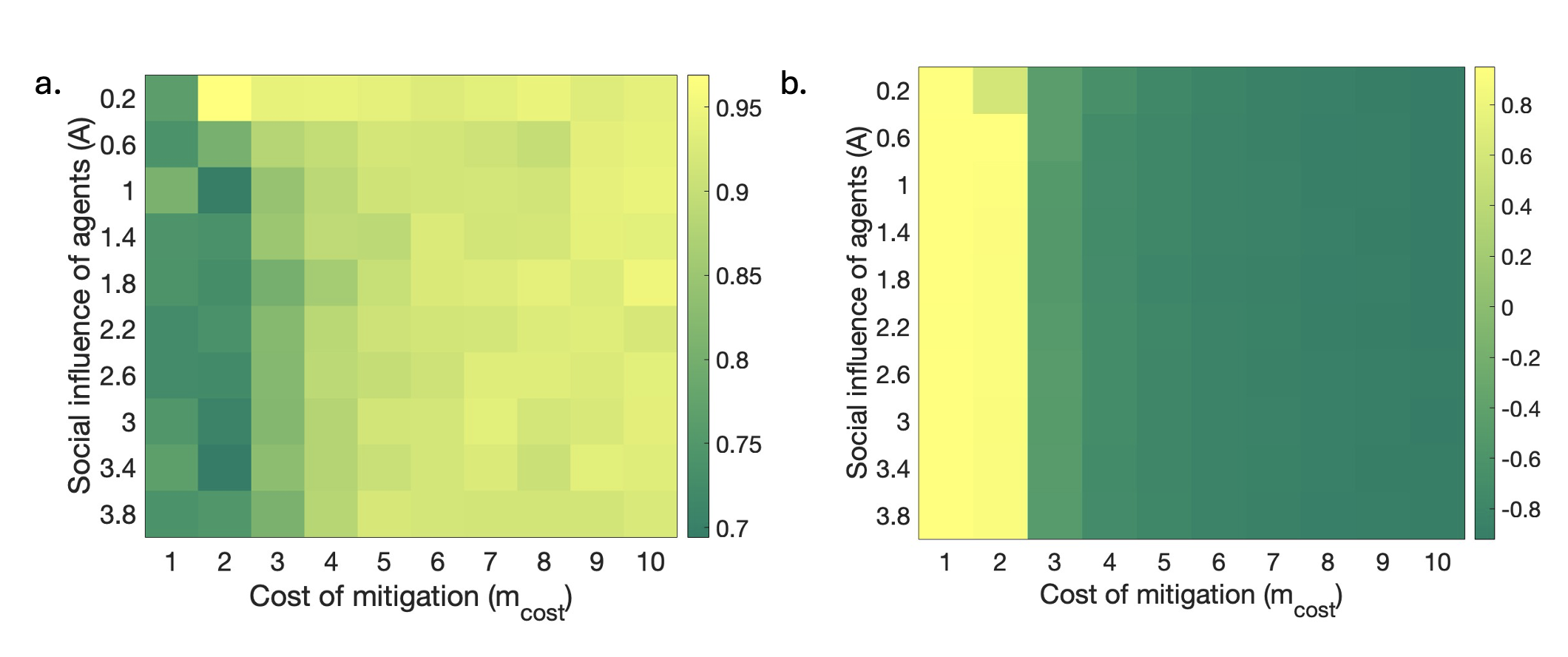}
\caption{\textbf{Social influence has minimal influence on average opinion with high mitigation costs: }Bivariate sensitivity of (a) bimodality coefficient, (b)average opinion when social influence($A$) is varied with the maximum response to temperature change $(R_{max})$}
\label{fig:B12}
\end{figure}
\begin{figure}[H]
\centering
  \includegraphics[width=18cm, height=7cm]{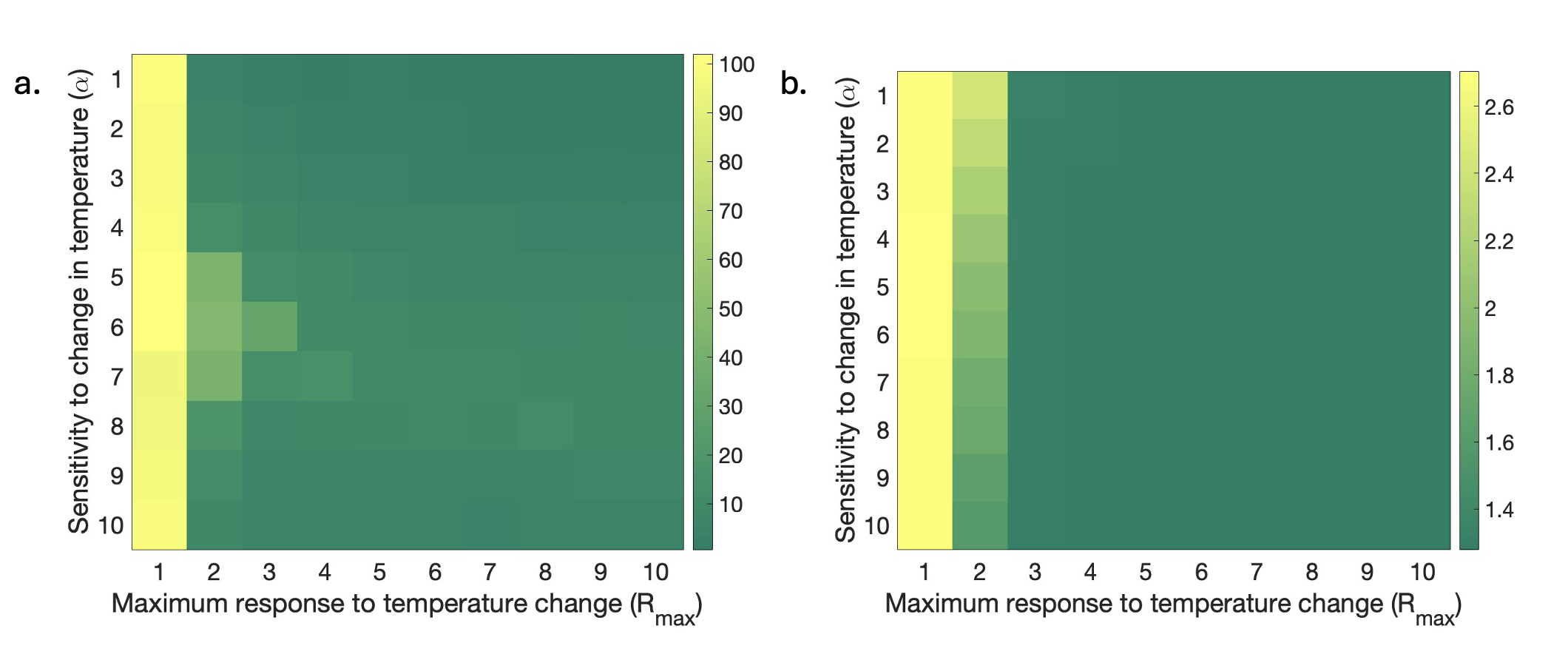}
\caption{\textbf{Sensitivity to temperature change has minimal influence on peak temperatures with high response: }Bivariate sensitivity of (a) emissions (GTCO\textsubscript{2}yr\textsuperscript{-1}), (b) peak temperature ($^{\circ}$C) when sensitivity to the temperature change ($\alpha$) is varied with the maximum response to temperature change $(R_{max})$}
\label{fig:B13}
\end{figure}
\begin{figure}[H]
\centering
  \includegraphics[width=18cm, height=7cm]{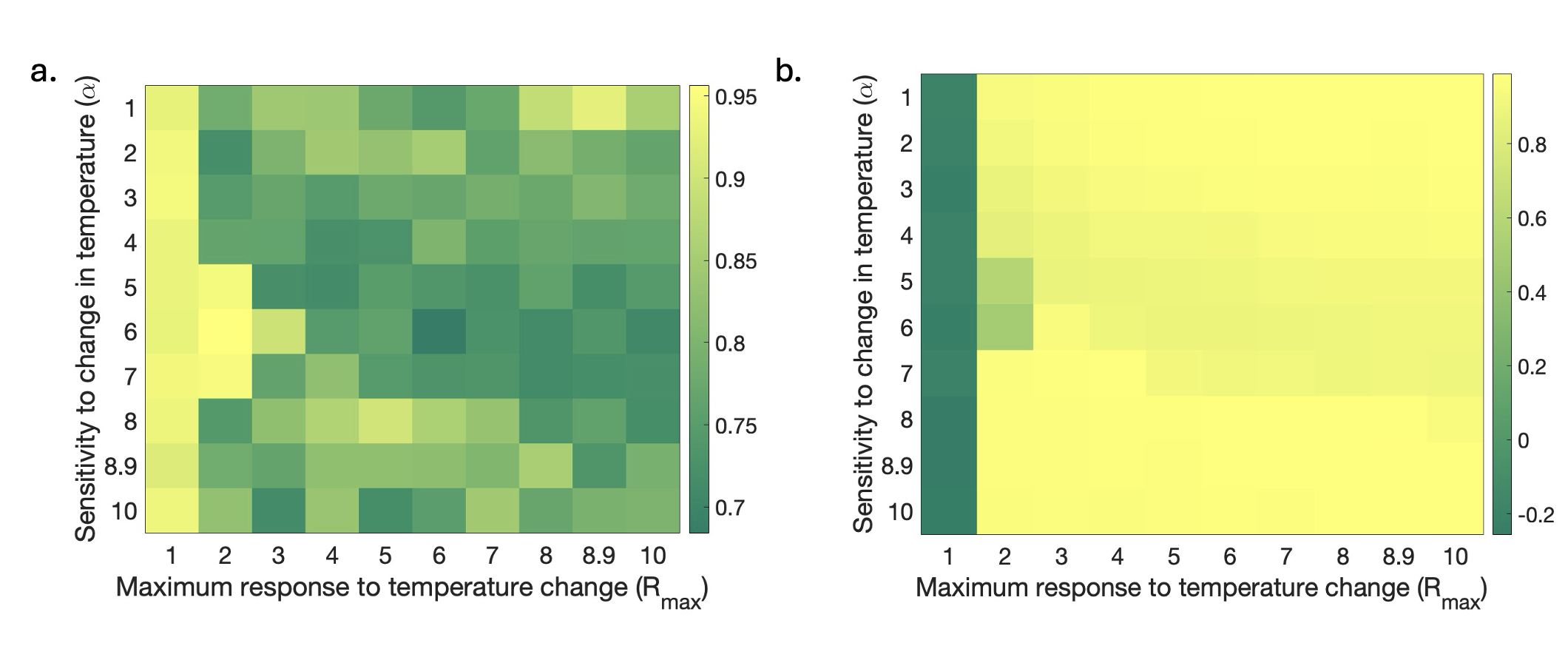}
\caption{\textbf{Sensitivity to temperature change has minimal influence on emissions with high response: }Bivariate sensitivity of (a) bimodality coefficient, (b) average opinion when sensitivity to the temperature change ($\alpha$) is varied with the maximum response to temperature change $(R_{max})$}
\label{fig:B14}
\end{figure}
\begin{figure}[H]
\centering
  \includegraphics[width=18cm, height=7cm]{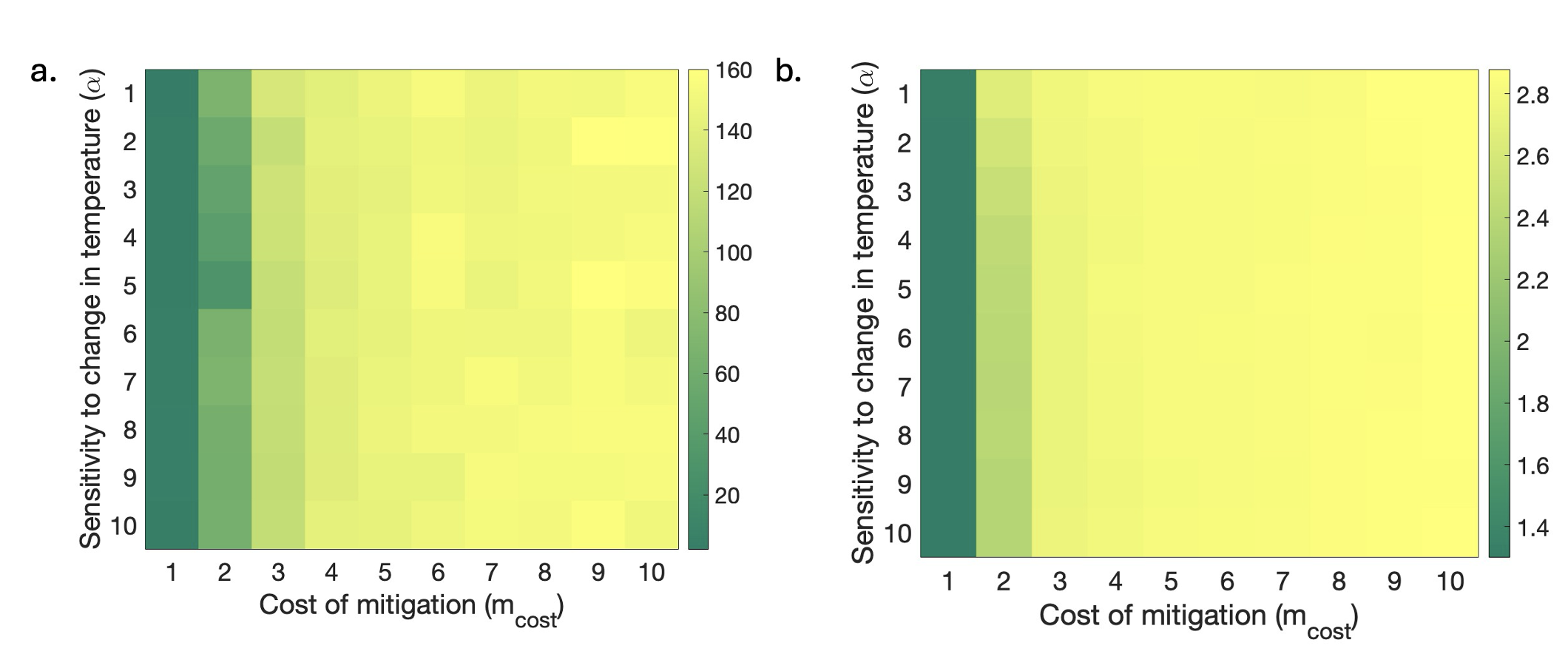}
\caption{\textbf{Sensitivity to temperature change has minimal influence on peak temperatures with high mitigation costs: }Bivariate sensitivity of (a) emissions (GTCO\textsubscript{2}yr\textsuperscript{-1}), (b) peak temperature ($^{\circ}$C) when sensitivity to the temperature change ($\alpha$) is varied with the cost of mitigation $(m_{cost})$}
\label{fig:B15}
\end{figure}
\begin{figure}[H]
\centering
  \includegraphics[width=18cm, height=7cm]{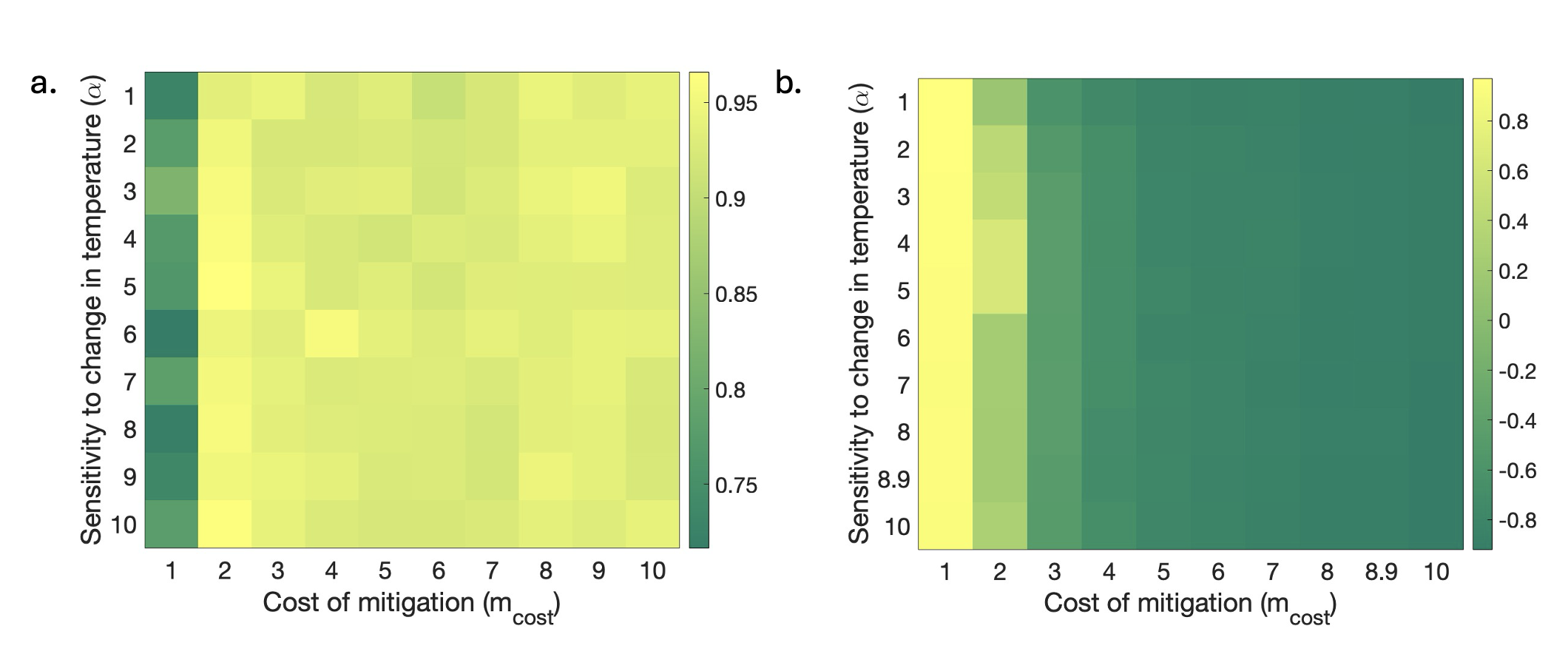}
\caption{\textbf{Sensitivity to temperature change has minimal influence on polarization with high mitigation costs: }Bivariate sensitivity of (a) emissions (GTCO\textsubscript{2}yr\textsuperscript{-1}), (b) peak temperature ($^{\circ}$C) when sensitivity to the temperature change ($\alpha$) is varied with the cost of mitigation $(m_{cost})$}
\label{fig:B16}
\end{figure}
\end{document}